\def\IR{{\mathbb R}} 
\def\oR{\overline {\IR}} 
\def\IS{{\mathbb S}} 
\def\IB{{\mathbb B}}

\def\IC{\mathbb C} 
\def\ID{{\mathbb D}}
\def\oC{\hat{\IC}}

\newcommand{\abs}[1]{\lvert1\rvert}
\newcommand{\Abs}[1]{\Big\lvert1\Big\rvert}
\newcommand{\norm}[1]{\lVert1\rVert}

\documentclass[12pt]{article} 

\usepackage[psamsfonts]{amssymb} 
\usepackage{amsfonts,amsmath}

\usepackage{epsfig,multicol}

\newtheorem{theorem}{Theorem}
\newtheorem{lemma}{Lemma}

\title{Frederick W.  Gehring  {\small (7 August 1925--29 May 2012) } \\{\small  Elected to NAS 1989}  }
\author{Gaven J Martin } \date{}
\begin{document}

\maketitle 
 \begin{center}
\scalebox{0.53}{\includegraphics*{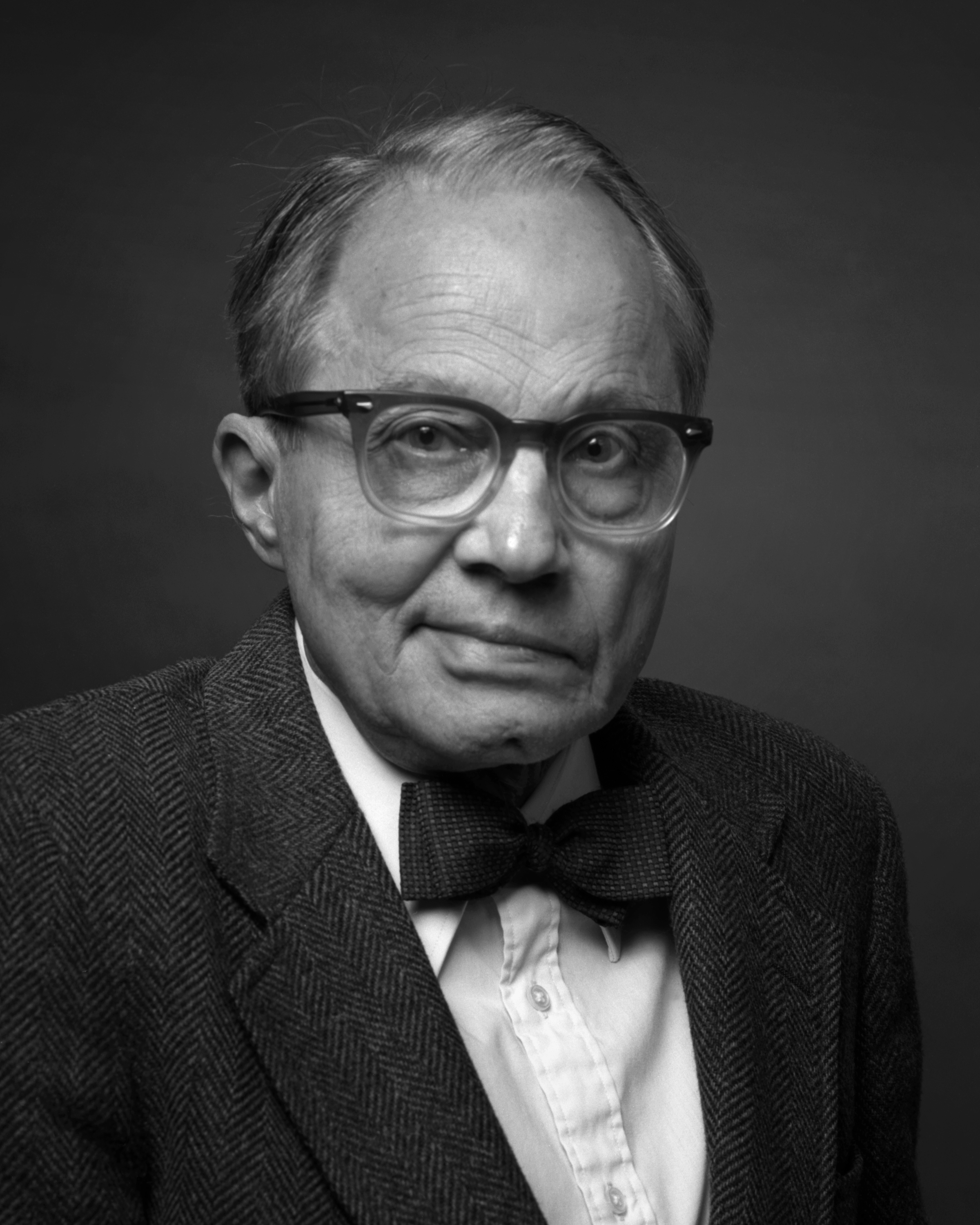}}\\
Fred Gehring
\end{center}
 \newpage

\begin{abstract}  Frederick William Gehring was a hugely influential mathematician who spent most of his career at the University of Michigan -- appointed in 1955 and as the T H Hildebrandt Distinguished University Professor from 1987.   Gehring's major research contributions were to Geometric Function Theory,  particularly in higher dimensions $\IR^n$,  $n\geq 3$.   This field he developed in close coordination with colleagues,  primarily in Finland,  over three decades 1960 -- 1990.  Gehring's seminal work drove this field forward initiating important connections with geometry and nonlinear partial differential equations,  while addressing and solving major problems.  During his career Gehring received many honours from the international mathematical community.  He was invited three times to address the International Congress of Mathematicians, at Moscow in 1966, at Vancouver in 1974, and at Berkeley (a plenary lecture) in 1986.  He was awarded honorary degrees from the University of Helsinki (1979), the University of Jyv\"askyl\"a (1990), and the Norwegian University of Science and Technology (1997).  In 1989, he was elected to the American Academy of Arts and Sciences and the National Academy of Sciences.  Other honors include an Alexander von Humboldt Stiftung (1981--84), Commander of the Order of Finland's White Rose (Commander class, Finland's highest scientific honour for foreigners), and a Lars Onsager Professorship at the University of Trondheim (1995--96).  He served for 19 years on various Committees of the American Mathematical Society.  In 2006 he received the American Mathematical Society's Steele Prize for Lifetime Achievement. He served three terms as Chairman of the Mathematics Department at the University of Michigan and played a leading role in shaping that department through the latter part of the 20th century.
\end{abstract}

\noindent {\bf Life.}  Frederick William Gehring was born  in Ann Arbor, Michigan on August 7, 1925.  His family was of German origin; his great-grandfather Karl Ernst Gehring (1829--1893) had emigrated from Germany in 1847 and settled in Cleveland, Ohio, where he founded the Gehring Brewery.  Gehring's grandfather 
Frederick William Gehring (1859--1925) was treasurer of the brewery and co-founder of a bank in Cleveland.  Gehring's father Carl Ernest Gehring (1897--1966) loved music and was an amateur composer.  He came to Ann Arbor to study engineering at the University of Michigan (UM), but soon switched to journalism and later worked for the {\it Ann Arbor News} as state news editor and music critic.  Gehring's mother Hester Reed Gehring (1898--1972) was the daughter of John Oren Reed (1856--1916), a physics professor at UM who later became Dean of the College.  She and Carl met as undergraduates at UM.   After the birth of their three children, she went on to complete a PhD in German and served as a foreign language examiner for the UM Graduate School.  

Gehring grew up in Ann Arbor and graduated from University High School in 1943.  He was then admitted to MIT to study physics or engineering, but he chose instead to enlist in the US Navy V-12 program, not knowing where he would be sent.  By coincidence, the Navy sent him to Ann Arbor for a special program in electrical engineering at UM.  He later graduated with a double major in electrical engineering and mathematics.  By that time the war was over, but the Navy sent him to sea for 4 months. 
      
      Upon his return, he re-enrolled at UM and decided, at the suggestion of Ruel Churchill, to concentrate on mathematics.  After receiving an MA degree from UM in 1949, he went to the University of  Cambridge  on a Fulbright Scholarship earning a PhD in mathematics in 1952, taking courses from such famous analysts such as J E Littlewood and A S Besicovitch while writing a thesis under the direction of J C  Burkill.   
      
      At Cambridge, Gehring met Lois Bigger, who had come from Iowa, also on a Fulbright Scholarship, to continue her study of microbiology. She also received her Ph.D. from Cambridge in 1952.  Gehring went to Harvard as a Benjamin Peirce Instructor, Bigger to Yale on a Research Fellowship.  They were married on August 29, 1953 in Bigger's hometown of Mt Vernon, Iowa.  
      
       \begin{center}
\scalebox{0.53}{\includegraphics*{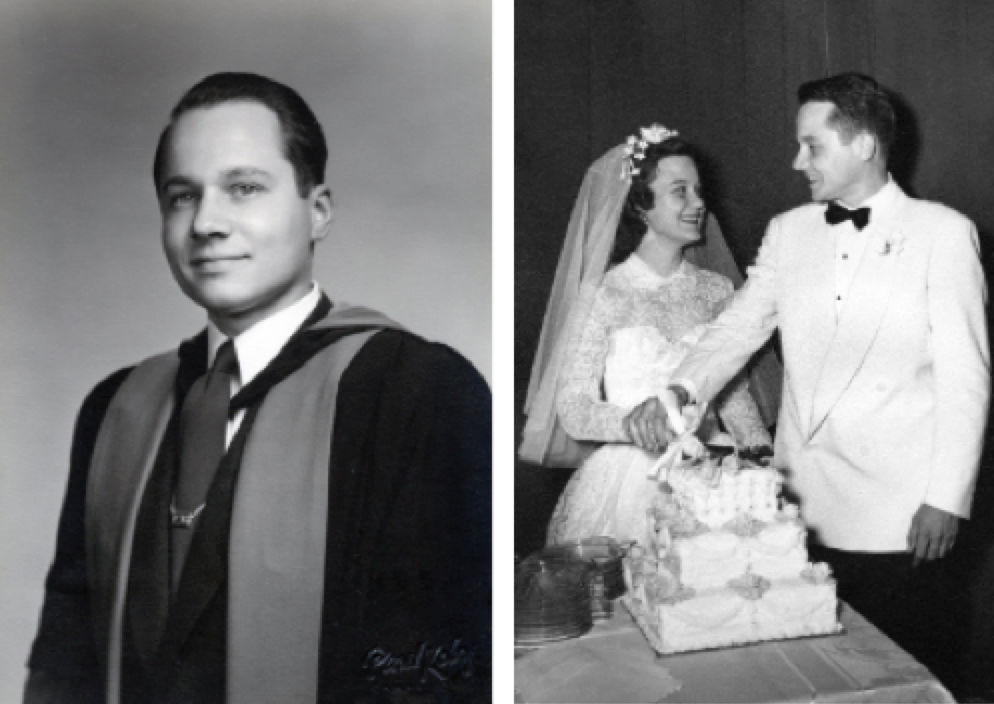}}\\
Gehring at graduation and wedding to Lois Bigger
\end{center}
      
      \bigskip
      
Gehring's acquaintance with Lars Ahlfors during his time at Harvard drew him closer to the area of complex analysis.  When the position at Harvard ended in 1955, T H Hildebrandt hired Gehring back to Ann Arbor, where he remained for the rest of his life.   A turning point in Gehring's mathematical career occurred in 1958.  Hoping to spend a year abroad, Gehring applied for Fulbright, Guggenheim, and NSF fellowships and was awarded all three.  His colleague Jack Lohwater had spent the year 1956--57 in Helsinki, where he saw that Olli Lehto and Kalle Virtanen, among others, were working on the theory of normal functions.  Encouraged also by Ahlfors, Gehring decided to go to Helsinki.  Lohwater put him in touch with Lehto, and Gehring arranged to spend the year 1958--59 there.  
      
Lehto \cite{lehto} tells the story of Gehring's year in Helsinki.  When Gehring arrived, he learned to his dismay that his Finnish hosts were no longer interested in meromorphic functions, but had been working on quasiconformal mappings.  As the story goes, Gehring said ``Fine,  I like quasiconformal mappings -- what are they?''  He was told he would learn, as they were about to start a seminar on the topic, in Finnish!  Gehring responded to the challenge and learned to speak Finnish that year.  Fred and Lois Gehring's first son Kalle was born in Helsinki that year, while their son Peter was born two years later in Ann Arbor.      
      
      In the year 1959--60,  Gehring moved to Z\"urich, where Albert Pfluger had also been working on quasiconformal mappings.  Stimulated by discussions with Pfluger and a paper by Charles Loewner, 
Gehring began to further develop the  higher-dimensional theory of quasiconformal mappings.  In those years abroad  Gehring  made professional contacts  which would strongly influence the entire course of his career.   When he returned to Michigan in 1960, Gehring began training students in quasiconformal 
mappings, and his first student graduated from UM in 1963.    In the course of his career, Gehring  directed 29 PhD students,  most of whom have had active careers in teaching and research at academic institutions. 
        
        An important aspect of Gehring's professional work was his extensive editorial service.  At various times in his career, he worked on the Editorial Boards of 9 different research journals.   He also served as Editor of book series for Van Nostrand (1963--70), North Holland (1970--94), and most famously for Springer-Verlag (1974--2003) with the Undergraduate Texts in Mathematics series. 
        
       At the University of Michigan, Gehring was promoted to Professor in 1962.  He served three terms as Chairman of the Mathematics Department, in 1973--75, 1977--80, and 1981--84.  He was named to a  collegiate chair in 1984, and became the TH Hildebrandt Distinguished University Professor in 1987.   The University honoured him with a Distinguished Faculty Achievement Award in 1981, the Henry Russel Lectureship in 1990, and a Sokol Faculty Award in 1994.  An international conference on ``Quasiconformal Mappings and Analysis" was held in Ann Arbor in August 1995 on the occasion of his 70th birthday.  He retired in 1996.
       
       Gehring was invited three times to address the International Congress of Mathematicians, at Moscow in 1966, at Vancouver in 1974, and at Berkeley (a plenary lecture) in 1986.  He was awarded honorary degrees from the University of Helsinki (1979), the University of Jyv\"askyl\"a (1990), the Norwegian University of Science and Technology (1997) and an Sc.D from Cambridge in 1976.  In 1989, he was elected to the American Academy of Arts and Sciences and the National Academy of Sciences (USA).  Other honours include an Alexander von Humboldt Stiftung (1981--84), Commander of the Order of Finland's White Rose (1988--2012), and a Lars Onsager Professorship at the University of Trondheim (1995--96).  He served for 19 years on various Committees of the American Mathematical Society, including 3 years on the Council, 4 years on the Executive Committee, and 10 years on the Board of Trustees.  He served on most of the AMS' important committees.  
              
       In 2006, the AMS honored Gehring with a Steele Prize for Lifetime Achievement.  The citation ({\it Notices Amer. Math. Soc.} {\bf 53} (2006), 468--469) says in part, {\em ``Largely because of Gehring's work, the theory of quasiconformal mappings has influenced many other parts of mathematics, including complex dynamics, function theory, partial differential equations, and topology.  Higher dimensional quasiconformality is an essential ingredient of the Mostow rigidity theorem and of recent work of Donaldson and Sullivan on gauge theory and four-manifolds\dots  Gehring's mathematics is characterized by its elegance and simplicity and by its emphasis on deceptively elementary questions which later become surprisingly significant."}
       
       \medskip
       
      Fred Gehring died in Ann Arbor on May 29, 2012 at age 86, after a long illness.   

\bigskip

\noindent {\bf Mathematics.}  In his lifetime Fred Gehring wrote around $130$ mathematical papers,  starting with his first  article in 1951. Gehring had wide and varied interests,  real and complex analysis,  geometric function theory, partial differential equations (pde), discrete groups and hyperbolic geometry,  and geometric group theory just to name a few.  He collaborated widely,  with nearly 40 different coauthors.  As noted above,  Gehring's PhD thesis at Cambridge was ostensibly under the supervision of J C Burkill,  but   he certainly felt that J E Littlewood and A S Besicovitch were equally involved in his research directions and mentorship. Indeed, he thanks Besicovitch for the problem which gave rise to his first published paper.

Any brief survey on Gehring's mathematical achievements  is going to miss a lot of important work.   There are obviously important contributions where he solved problems that had been open for decades,  but there are other papers which introduced key concepts and initiated other work by subsequent authors.  This of course does not even begin to address the many ideas and mathematically fruitful routes which he freely gave away to others,  and in particular his many students.    Gehring's many and varied projects were often running simultaneously so it is impossible to  order his work in time.  There are threads that run through his entire career and the study of quasidisks is one-such.  Gehring's paper {\em Quasiconformal mappings in Euclidean spaces} in the Handbook of complex analysis, {\bf 2}, 1--29, 2005 gives a fairly concise summary of his perspective of the main results achieved in the theory of quasiconformal mappings since the 1930's and his broad vision of how they connect with other areas of mathematics.  Gehring was certainly leading many of these developments.  

Gehring  was a meticulous writer,  carefully crafting the statement of a theorem so as to reflect the contents of the proof.  He worked tirelessly to ensure that an outline and an estimate became a sharp proof and he never used a result unless he understood its proof.  Thus each paper was typically accompanied by an ``Idiots Guide'' carefully explaining all the details,  many ran to several hundred pages to accompany a twelve page paper.   

\medskip

\noindent {\bf Early years.} Gehring's first paper is {\em Images of convergent sequences in sets},  published in  J. London Math. Soc. 26, (1951),  249--256.   It builds on earlier work of  H. Hadwiger (well known for his work on convex bodies) and from  H. Kestelman both in 1947.    Hadwiger showed that if a set $E\subset \IR^n $ has positive interior measure, then it contains a self similar image of any finite subset $S\subset \IR^n$.  Gehring considered this for infinite sets $S$, and in particular when  $S$ is a convergent sequence - somewhat closer to Kestelman's problem.   Gehring's first significant paper was that he wrote with Lehto on his first visit to Finland.  That paper is {\em On the total differentiability of functions of a complex variable}, appearing in
Ann. Acad. Sci. Fenn. Ser. A I, 1959 9 pp.  The paper is still a classic and is often used by  researchers,  see  \cite[\S3.3]{AIM}.  The result is false in higher dimensions and it is the interplay between topology and analysis in the proof that was one of Gehring's favourite things. The main theorem of the paper is the following:
 \begin{theorem}
Let  $f:\Omega\to \IC$ be a continuous open mapping of a planar domain $\Omega$. Then $f$ is differentiable almost everywhere in $\Omega$ if and only if $f$ has finite first
partials almost everywhere.
\end{theorem}
 \begin{center}
\scalebox{0.75}{\includegraphics*{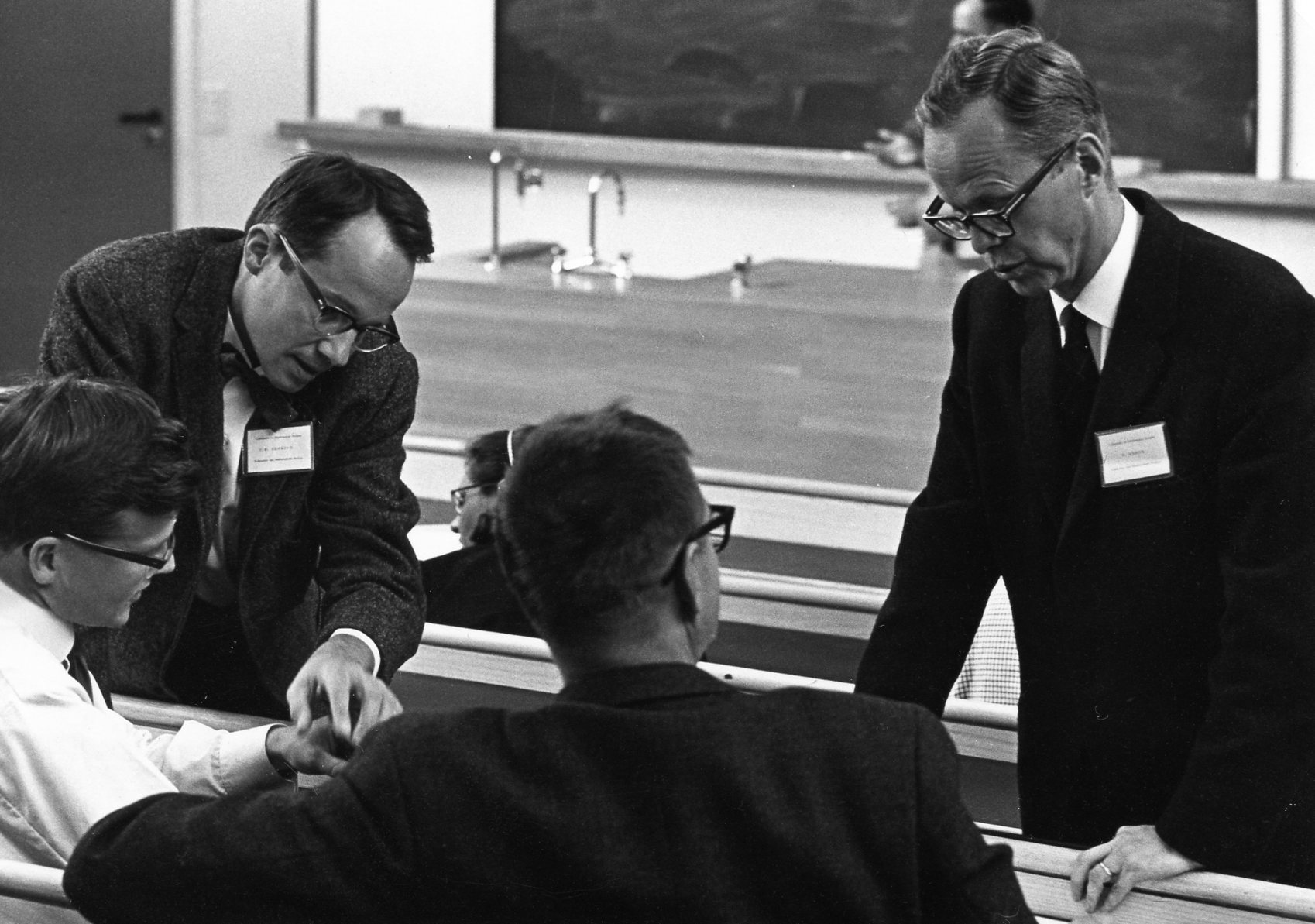}}\\
Gehring, with from left,  Jussi V\"as\"al\"a,  Lalevi Suominen and Olli Lehto 
\end{center}

For  planar homeomorphisms this result had been earlier established by D. Menchoff in 1931.  The proof that Gehring and Lehto give is a model of exposition and has not changed much over the years - it is almost basic measure theory,  but with a wonderful insight.  The topology comes in as open mappings satisfy a maximum principle.  As a consequence of this theorem every discrete open Sobolev mapping  $f \in W^{1,1}_{loc}(\Omega,\IC)$ is differentiable almost everywhere.  Using this result one can connect the volume derivatives and  the pointwise Jacobian of a mapping and thereby obtain a first versions of the change-of-variables formula and further  observe that for discrete open mappings the Jacobian cannot change sign on the domain of definition,  so analytic degree is defined.  Clearly this result is central,  underpinning a lot of planar geometric function theory.

\medskip

\noindent  {\bf Higher dimensional quasiconformal mappings.}  From early in 1960 Gehring began the task of laying down some of the (now) basic tools of quasiconformal mappings in higher dimensions.  There were others working closely in this area,  most notably Jussi V\"ais\"al\"a (in Finland) and Yu. G. Reshetnyak (in Russia).  Gehring wrote quite a long review of V\"ais\"al\"a's papers, {\em On quasiconformal mappings in space} and  {\em On quasiconformal mappings of a ball} (see \cite{V1,V2}).   The major tool in V\"ais\"al\"a's work was the  extremal length method - that is via the geometric definition of quasiconformality.  Gehring and V\"ais\"ala had already written together on this method in {\em On the geometric definition for quasiconformal mappings},  Comment. Math. Helv.,  36, (1961) 19--32.  Of course  these ideas for planar quasiconformal mappings go  back  to Lars V. Ahlfors \cite{A}, Pfluger \cite{P} and even before to Ahlfors' work with Arne Beurling in 1946 \cite{AB2}.    

Initially both Gehring and V\"ais\"al\"a were seeking to show the equivalence of the various natural generalisations to higher dimensions of the many equivalent  definitions for planar quasiconformal mappings - see two below.     They also sought extensions of other well known planar results to higher dimensions.  For instance V\"ais\"al\"a established the Caratheodory Theorem in this context in \cite{V2}.  Let us give two possible definitions for a quasiconformal mapping of a domain $\Omega\subset\IR^n$.  

Here then are two apparently independent definitions for a quasiconformal mapping. First the {\em analytic definition}.  Suppose $f$ is a homeomorphism which lies in the Sobolev space  $W^{1,n}_{loc}(\Omega,\IR^n)$ of functions whose first derivatives are integrable with exponent $n={\rm dimension}$.  Then $f$ is {\em $K$-quasiconformal}, $1\leq K<\infty$,  if
\begin{equation}\label{1st}
|Df(x)|^n \leq K \; J(x,f), \hskip10pt \mbox{for almost all $x\in \Omega$}.
\end{equation}
Here $Df$ is the $n\times n$ differential of $f$,  $J(x,f)$ is its determinant (Jacobian) and $|\cdot|$ is the operator norm.  

Next a {\em geometric definition}.  Let $\Gamma$ be a family of curves in $\Omega$.  For instance the family of all curves connecting two sets $E,F\subset \Omega$.  A non-negative Borel integrable function $\rho:\Omega\to\IR_{\geq 0}$ is admissible for $\Gamma$ if for each $\gamma\in \Gamma$ we have
\[ \int_\gamma \rho(s) \; ds \geq 1 \]
The modulus of $\Gamma$ is
\[ M(\Gamma) = \inf \Big\{ \int_\Omega \rho^n(x) \; dx : \mbox{$\rho$ is admissible for $\Gamma$} \Big\} \]
We say that a homeomorphism $f:\Omega\to f(\Omega)$ is {\em $\tilde{K}$-quasiconformal}, $1\leq \tilde{K}<\infty$,  if for every curve family $\Gamma\subset\Omega$ 
\begin{equation}\label{2nd}
\frac{1}{\tilde{K}}\; M(\Gamma) \leq M(f\Gamma) \leq \tilde{K}\; M(\Gamma)
\end{equation}
Here $f\Gamma=\{f\circ \gamma:\gamma\in \Gamma\}$ is a curve family in $f(\Omega)$.

Hadamard's inequality gives $|Df(x)|^n \geq   J(x,f)$ and so equality holds in (\ref{1st})  for $K=1$ and if also the dimension $n=2$,  a little calculation reveals the Cauchy-Riemann equations.  Further,  and relevant to our later discussion, from the Looman-Menchoff theorem we see the weaker hypothesis that $f\in W^{1,1}_{loc}(\Omega,\IR^2)$ would suffice to conclude $f$ is conformal.  For (\ref{2nd}), when $\tilde{K}=1$, the argument is a little longer to deduce that $f$ is conformal but the relationship to the length-area inequality for conformal mappings is clear.  Next,  in dimension $n\geq 3$  the constants $K$ and $\tilde{K}$ may differ for the same mapping, but both are simultaneously equal to $1$ (not hard) and  simultaneously finite or  infinite (a hard theorem).  Thus the space of quasiconformal mappings is well defined,  but when specifying constants one must refer to a particular definition.

\medskip

Gehring's two major papers here both appeared in the Transactions of the American Math. Soc,  {\em Symmetrization of rings in space}, TAMS,  101,  (1961) 499 -- 519 and  {\em Rings and quasiconformal mappings in space}, 
TAMS, 103, (1962), 353--393.  There was an earlier announcement of the main results in the Proceedings of the National Academy,  {\bf  47}, (1961), 98--105.    Gehring starts with results on the spherical symmetrisation of rings in space.  This is quite difficult,  but motivated by earlier two dimensional results of G  Bol, T  Carleman and  G  P\'olya and G  Szeg\"o among others.  This symmetrisation basically  identifies extremal configurations for rings (read conformal invariants) in terms of geometric information such as diameters of components and distance between them - generalising the ``length-area'' method from complex analysis.  From this one can give estimates on arbitrary rings.  This is a key fact and leads to direct proofs of such things as the H\"older continuity of higher dimensional quasiconformal mappings previously established by E D Callender \cite{Call}, following arguments of C B Morrey and L  Nirenberg).  Naturally,  this modulus of continuity will give compactness via the Arzela-Ascoli theorem.  Next,  Gehring proved the existence and uniqueness of an extremal function realising the Loewner capacity of a ring (doubly connected domain) in space.  This extremal function $u$ solves the nonlinear pde ${\rm div}(|\nabla u|\nabla u)=0$.   Gehring then used these results to establish equivalences between the analytic definition  and the geometric definition.  These proofs establish deep connections between the use of the Rademacher-Stepanoff theorem and the absolute continuity and geometric properties of mappings.

The most recognisable thing that Gehring proved using these two papers is the Liouville theorem.  In 1850,  Joseph Liouville added a 
short note to a new edition of Gaspard Monge's classic work 
{\em Application de l'Analyse \`a la G\'eometrie\/}, whose publication 
Liouville was overseeing.
The note was prompted by a series of three letters that Liouville had 
received in  1845/6 from the  British physicist 
William Thomson (better known  as Lord Kelvin) who had studied in Paris under 
Liouville's  in the mid-1840s.   In his letters, Thomson asked Liouville a number of questions concerning
inversions in spheres, questions that had arisen in conjunction with 
Thomson's research in electrostatics, in particular, with the so-called 
principle of electrical images 
(note that  reflection in the unit sphere $\IS^2$ of $\IR^3$
is often referred to in physics as the ``Kelvin transform.'')   In Liouville's time conformal mappings were certainly understood  to be many times differentiable and following his motivation for writing the article, Liouville framed
his discussion in the language of differential forms so his original formulation  bears little resemblance
to the modern-day  theorem. However Liouville's title,  ``Extension au cas de trois dimensions de la question du 
trac\'e g\'eographique''  gives no hint whatsoever as to the results.  It was only later that Liouville published his theorem in a form 
approximating the usual statement which we now discuss.  

In higher dimensions the group of all M\"obius or Conformal transformations of
$\oR^n$  consists of all finite compositions of reflections in spheres and hyperplanes.  Thus it contains the similarities $x\mapsto \lambda Ox+b$,  $\lambda >0$,  $O\in O(n)$,  the orthogonal group,  and $b\in \IR^n$ and the inversions $x\mapsto \frac{x-a}{|x-a|^2}$.  It is easy to see that these mappings provide 
examples of conformal transformations - their derivative at any point is a scalar multiple of an orthogonal transformation and therefore the map infinitesimally preserves angles (i.e. is conformal).  These mappings are  all infinitely differentiable and the space of all M\"obius transformations of $\oR^n$ forms a {\em finite dimensional} Lie group.
 
Liouville proved in 1850 that if $f:\Omega\to\IR^3$ is a 3 times continuously differentiable conformal mapping,  then there is a M\"obius transformation $\Phi:\oR^n\to\oR^n$ such that $\Phi|\Omega=f$.  This is a very surprising rigidity theorem as in two dimensions the space of smooth conformal mappings is infinite dimensional and there is no local-to-global extension theorem.  One significant corollary of Liouville's theorem is that
the only subdomains $\Omega$ in $\IR^n$ with $n\geq 3$ that are conformally 
equivalent to the unit ball $\IB^n$ are Euclidean balls and half-spaces.  This stands in stark contrast to the famous Riemann mapping theorem, announced in 1851 
a year after Liouville's note was published:
any simply connected proper subdomain $\Omega\subset \IC$ of the complex 
plane is conformally equivalent to the unit disk $\ID$.

Certainly conformal mappings are $1$-quasiconformal mappings - that is conformal invariants such as moduli are preserved.  What of the converse.  Certainly this is implied when Liouville's assumptions on smoothness hold.  Gehring had by now developed the tools to identify extremal rings and their configurations and an understanding of how modulus behaved under continuous deformations and how this gave regularity for quasiconformal mappings.  It was a short step (but the one he had been aiming at) to prove the following theorem.

\begin{theorem}\label{LT} Let $\Omega\subset \IR^n$, $n\geq 3$ and let $f:\Omega\to f(\Omega)$ be a $1$-quasiconformal mapping. Then $f$ is the restriction to $\Omega$ of a M\"obius transformation of $\oR^n$.  
\end{theorem}
 
There are now refinements of this result,  see \cite{I1,I2} and the curious discrepancy between what is known in even and odd dimensions remains one of the most interesting problems in the higher dimensional theory.
 
 \medskip

Lars  Ahlfors' reviews of these papers of Gehring begins {\em ``This is an announcement of important results in the theory of 3-dimensional quasiconformal mappings''} and ends  ``{\em In a final section he [Gehring] proves that a $1$-quasiconformal mapping is a M\"obius transformation. This is Liouville's theorem without any regularity assumptions, a remarkable achievement.''}

\medskip

This last note was professionally important to Gehring as certainly Lars Ahlfors was a personal and mathematical hero to him,  perhaps above all others.

 \begin{center}
\scalebox{0.75}{\includegraphics*{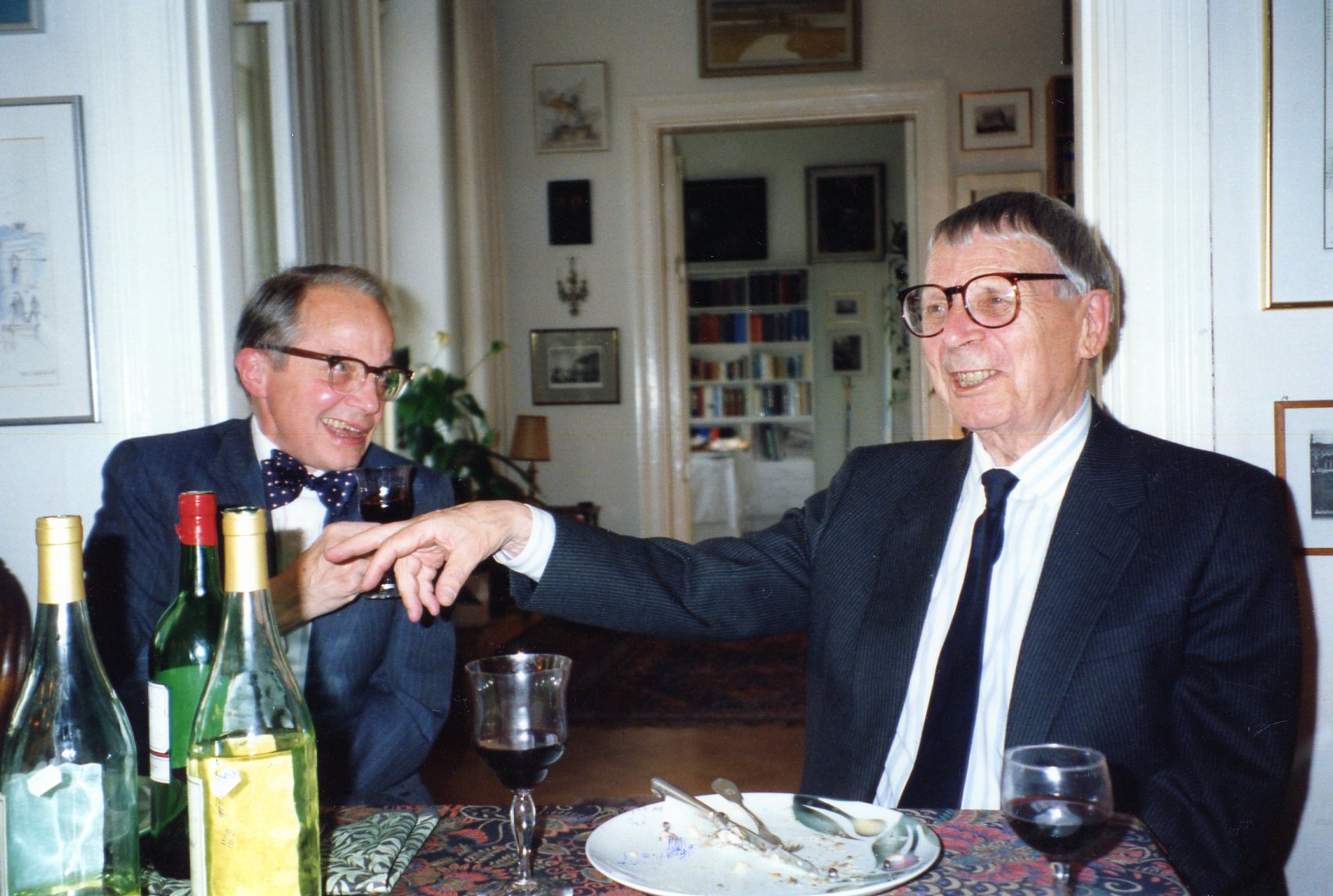}}\\
Gehring and Lars V Ahlfors \end{center}
 
 \medskip
 
Soon after these papers were published Gehring and  V\"ais\"al\"a worked together on the paper {\em The coefficients of quasiconformality of domains in space}, in  
Acta Math., 114 (1965) 1--70.   The problems they consider are typically of the following form:  Suppose $\Omega\subset\IR^n$ is a domain.  Is there a quasiconformal homeomorphism $f:\IB^n\to\Omega$ ? and if so,  bound its distortion.  They had worked together on these  sort of problems before,  but this paper produces the most definitive results (and led to a number of students  extending the results as well,  G Anderson, M.K. Vamanamurthy and K. Hag wrote on this subject,  see \cite{And,Hag1,Hag2}).   As noted, a problem here is that there are different measures of distortion depending on the definition of quasiconformality used.  Gehring and V\"ais\"al\"a carefully define the different distortion functions (which are all quantitatively,  but not functionally related). Precise results depend on the distortion chosen. 

In two dimensions all distortion functions coincide and a considerable literature on extremal mappings and  various existence criteria due to the important  relationship with Teichm\"uller theory.  In higher dimensions these issues are still not resolved (not even close actually).  
The  calculation of the coefficients of the dihedral wedge, cylinder, and cone by Gehring and V\"ais\"al\"a were the first examples of the calculation of such distortion quantities in the literature.    Perhaps the most well known results of the paper,  and certainly one of the most useful,  is the following.
\begin{theorem} Suppose that $0<a<b$, that $D$ is a domain in $\IR^3$, and that $C(D)\cap\{|x|<b\}$ has at least two components which meet $\{|x|=a\}$. Then the {\em inner distortion},  $K_I(D)\geq A\log(b/a)$. $A\geq 0.129$ is an absolute constant.
\end{theorem}

An interesting case they consider is the  infinite cylinder $D=\{x=(r,\theta,x_3)\colon 0\leq r<1\}$.  They prove the {\em outer distortion} $K_0(D)=(q/2)^{1/2}$, where $q$ is the EllipticK functional value
\[ q=\int_{0}^{\pi/2}\ \frac{1}{\sqrt{\sin u}} \; du \]
 The extremal map is obtained as follows. $D$ is mapped onto the half-space $D'=\{x_3>0\}$ by $f_1(r,\theta,x_3)=(t,\theta,\varphi)$, where $(t,\theta,\varphi)$ are spherical coordinates, $$ r=\left(\frac 1{q}\int_0^\varphi\ (\sin u)^{-1/2}du\right)^2,\quad x_3=\frac 2{q}\log t. $$ Then $D'$ is mapped onto  the unit $3$-ball by a M\"obius transformation.  Of course the hard part is explicitly computing certain conformal invariants in the domain and range to show this estimate is sharp.
 
 Next they prove the following:

\begin{theorem}
If $\partial D$ contains an inward directed spire or an outward directed ridge, then $K(D)=\infty$.   For each $\varepsilon>0$ there exists a domain whose boundary contains an outward directed spire and whose coefficients are within $\varepsilon$ of the corresponding coefficients of an infinite cylinder.
\end{theorem}  
They use this result to construct a quasiconformal ball whose set of non-accessible boundary points has positive 3-dimensional measure.  

\medskip

There is important recent work related to this discussion,  but obtaining nice geometric criteria on a topological ball to guarantee it is a quasiball is probably never going to happen.    However many necessary criteria for a domain to be quasiconformally equivalent to a ball have found wide applications.  The notion of uniform and John domains,  quantitative versions of local connectivity and so forth have become important ideas in the theory of analysis on metric spaces as laid down by J. Heinonen and P. Koskela \cite{HK},  and elsewhere.

\scalebox{0.7}{\includegraphics*[viewport=-100 10 300 240]{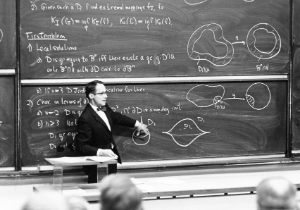}}
 \begin{center}
 Gehring lecturing on quasiballs
 \end{center}
\medskip

\noindent {\bf Area distortion.}  Until Kari Astala's solution in 1994,  \cite{Astala},  the most famous open problem about planar quasiconformal mappings was the area distortion conjecture raised in Gehring's paper with Edgar Reich from 1966,  {\em Area distortion under quasiconformal mappings},   Ann. Acad. Sci. Fenn. Ser. A I (1966) 15 pp.  Suppose that $f(\ID)=\ID$ and $f(0)=0$.  Given $E\subset \ID$,  denote it's measure by $|E|$.  It follows from Bojaski's improved regularity result for planar quasiconformal mappings (see the next section,  and  also \cite{Bo}) that for each $K\geq 1$ there are numbers $a_K$ and $b_K$, such that if $f$ is $K$-quasiconformal, then
\begin{equation}\label{1}
 \frac{ |f(E)| }{\pi}\leq  b_K \Big(\frac{|E|}{\pi}\Big)^{a_K}, \hskip15pt\mbox{ for each measurable set $E\subset \ID$}
  \end{equation}
Gehring and Reich prove $a(K)=K^{-a}$ for some  $1\leq a\leq 40$, and that $b(K)=1+O(K-1)$ as $K\rightarrow 1$. They deduce this as a consequence of the inequality for the Beurling transform of characteristic functions
\begin{equation}\label{2} \iint_\ID |{\cal S}\chi_E(z)|\,dxdy\leq a \, |E| \log(\pi/|E|)+b\, |E| \end{equation} 
for some $1\leq a\leq 40$, and $0\leq b\leq 2$ which they achieve using  the Calder\'on-Zygmund theory of singular integral operators.  They conjectured that $a=1$ in both the inequalities (\ref{1}) and (\ref{2}) and showed that the best bound in either (\ref{1}) or (\ref{2}) gives the best bound in the other.  Astala showed,  among other things, that this first area distortion conjecture is true.  But the other problem suggested by this paper concerns the $L^p$-norms of the Beurling transform,  and that still stands as one of the most challenging problems in the modern theory.

\medskip
 
\noindent {\bf Higher integrability.}  We now come to  the result that is undoubtedly Gehring's most well known contribution to mathematics.  Though the results have been improved,  many of the ideas are still used in mathematics today.  
 
In a remarkable paper in 1973,  actually a landmark in modern analysis and pde,   {\em The $L^{p}$-integrability of the partial derivatives of a quasiconformal mapping},  
Acta Math., 130 (1973), 265--277,  Gehring established that the Jacobian determinant of a
$K$-quasiconformal mapping is integrable above the natural exponent ($n={\rm dimension}$).  That is the usual assumption $f\in
W^{1,n}_{loc}(\Omega,\IR^n)$,   together with a bound as per (\ref{1st}) implies that $f\in
W^{1,n+\epsilon}_{loc}(\Omega,\IR^n)$ for some $\epsilon>0$ depending on $n$ and $K$ and for which Gehring gave explicit estimates.

While this result was already known in the plane due to the work of B  Bojarski
\cite{Bo}, and perhaps anticipated in higher dimensions, it is impossible to overstate how important this result 
has proven to be in the theory of quasiconformal mappings and more generally Sobolev spaces and in the
theory of non-linear pdes.  The techniques developed to solve this problem,  for instance the
well--known reverse H\"older inequalities, are still one of the main tools used in  non-linear potential theory, non-linear elasticity, pdes and harmonic analysis. In the introduction to his paper,  Gehring says  the results follow  {\em ``using a quite elementary proof''}. P  Caraman repeats this claim in his review, however time has shown Gehring's ideas here to be deep and highly innovative.  It is perhaps a mark of the clarity of Gehring's writing (and his own modesty) that one can follow the arguments easily and overlook how hard these ideas are to come by.

\medskip

We state the following version of Gehring's result as proved in \cite{IM} which also gives the result  for quasiregular mappings.
\begin{theorem}[Higher Integrability]  Let $f:\Omega\to\IR^n$ be a mapping of Sobolev class $W^{1,q}_{loc}(\Omega)$ satisfying the differential inequality
\begin{equation}\label{difineq}
|Df(x)|^n \leq K\; J(x,f).
\end{equation}  
Then there are $\epsilon_{K}{*},\epsilon_K>0$ such that if $q>n-\epsilon_{K}{*}$,  then $f\in W^{1,p}_{loc}(\Omega)$ for all $p<n+\epsilon_K$.
\end{theorem}
As an immediate corollary we have the following:
\begin{theorem}\label{GHI}  Let $f:\Omega\to\IR^n$ be a $K$ quasiconformal mapping.  Then there is $p_{_K}>n$ such that  $f\in W^{1,p_{_K}}_{loc}(\Omega)$.
\end{theorem}
The higher-dimensional integrability conjecture,  perhaps the most important outstanding problem in the area,    would assert that if $f$ satisfies (\ref{difineq}) and lies in $W^{1,q}_{loc}(\Omega)$ for some $q>nK/(K+1)$, then $f$ actually lies in the Sobolev space $W^{1,p}_{loc}(\Omega)$ for all $p<nK/(K-1)$.  In two-dimensions this conjecture was proven by Astala as mentioned earlier.  

A key step in Gehring's proof is   the following\index{reverse H\"older inequality}
\begin{lemma}[Reverse H\"older inequality]  Let $f:\Omega\subset \IR^n \to \Omega'$ be a $K$-quasiconformal mapping.  Then there is $p=p(n,K)>1$ and $C=C(n,K)$ such that
\begin{equation}\label{RHI}
\Big(\frac{1}{|Q|}\, \iint_Q \; J(x,f)^p \; dx \Big)^{1/p} \leq \frac{C}{|Q|}\, \iint_Q \; J(x,f) \; dx
\end{equation}
for all cubes $Q$ such that $2Q\subset\Omega$.
\end{lemma}  
This result is often referred to as the {\em Gehring Lemma}.  It is  worth reading T. Iwaniec's survey on this lemma \cite{IwaniecGL}. That paper contains many recent related results, with numerous generalizations and some of the far-reaching applications in mathematical analysis.

Next,  one deduces from the improved regularity bounds on the distortion of Hausdorff dimension ${\rm dim}_H(E)$ of a set $E\subset \IR^n$.
\begin{theorem}
If $f:\Omega \rightarrow \IR^n$ is $K$-quasiconformal,  then there is $a\geq 1$ such that for each $E\subset\IR^n$, 
\[ K^{-a}\;  \Big( \frac{1}{{\rm dim}_H(f(E))} -\frac{1}{n}  \Big) \leq \frac{1}{{\rm dim}_H(f(E))} -\frac{1}{n}  \leq K^a \Big( \frac{1}{{\rm dim}_H(f(E))} -\frac{1}{n}  \Big) \] 
\end{theorem}
This result shows in particular that sets of $0$ and $n$ dimensional Hausdorff measure are preserved,  proving an earlier conjecture Gehring made in his work with V\"ais\"al\"a.  The optimal conjecture would be that $a=1$ here.

\medskip

Gehring was always keen to give a lot of credit to E  De Giorgi's results \cite{DG} to support his proof.   While De Giorgi's work (and also that of J F Nash) was hugely important,  it was substantially earlier and well known in the pde community by that time.  Gehring found considerable insight and gave many new ideas to solve these problems.

\medskip

\noindent {\bf Quasiconformal Homogeneity.}  Another theme that continues to generate a lot of interest is Gehring's paper with former student Bruce Palka,  {\em Quasiconformally homogeneous domains},  which appears in  J. Analyse Math. 30 (1976), 172--199.  As   mentioned earlier,  one of Gehring's lifelong tasks was to find a characterisation of those domains $\Omega\subset\IR^n$ which are of the form $\Omega=f(\IB^n)$ where $f$ is a quasiconformal mapping and $\IB^n$ is the unit ball - {\em quasiballs} (quasidisks in two dimensions).  This paper was written ``{\em in the hope that this may eventually lead to a new characterisation  for domains quasiconformally equivalent to $\IB^n$}''. There appears to be no simple characterisation such as Ahlfors' criterion in higher dimensions.  So some other ideas need to be advanced. There are clear characterisations in terms of quasiconformal generalizations of the local collaring theorems of M Brown and others from geometric topology.  Here one assumes  that the boundary $\partial \Omega$ has a neighborhood $U$ such that $U\cap\Omega$ can be mapped quasiconformally into $\IB^n$ so that $\partial \Omega$ corresponds to $\partial B^n$ (Gehring first did this in 1967.  It directly implies that smoothly bounded Jordan domains are quasiballs in higher dimensions).  Another criterion,  unfortunately very hard to apply,  is  that of L.G. Lewis \cite{Lewis}: the corresponding $n$-dimensional Royden algebras $A_n(\Omega)$ and $A_n(\IB^n)$ should be isomorphic.

In this work mentioned above,  Gehring and Palka  study domains $\Omega$ which are homogeneous with respect to a quasiconformal family, a quasiconformal group, or a  conformal  family, where homogeneous with respect to a family $\Gamma$ means that for each pair $a,b\in D$, there exists an $f\in\Gamma$ such that $f(a)=b$. 

Gehring and Palka also prove  that if $D$ is homogeneous with respect to a quasiconformal family and has a $k$-tangent $(1\leq k\leq n-1)$ at some finite boundary point, then $D$ is quasiconformally equivalent to (a) $B^n$ if $k=n-1$ and to (b) $\overline \IR^n-\overline \IR^k$ if $k=1,\cdots,n-2$, $n>2$. They obtain also, in the particular case $n=2$, that $D$ is homogeneous with respect to a quasiconformal group if and only if $D$ can be mapped conformally onto $\overline \IR^2$, $\IR^2$, $\IR^2-\{0\}$ or $\ID$.  In the $n$-dimensional case they obtain similar results when $D$ is homogeneous with respect to a $1$-quasiconformal family.  There the results basically follow from Lie theory since the $1$-quasiconformal mappings form a group.  Actually it turns out that quasiconformal groups are Lie groups in all dimensions as well \cite{MLG}.

 They give examples to show there can be no easy characterization for $\Omega$ in $\IR^n$ quasiconformally equivalent to $B^n$ in terms of the tangential or smoothness properties of their boundaries.

Recent work on qc-homogeneity has connections with rigidity.  The question of the existence of a lower bound (strictly greater than one) on the quasiconformal homogeneity constant of quasiconformally homogeneous sets first appeared in this paper.  Recently F. Kwakkel and V. Markovic prove
that here exists a constant $\epsilon$ such that, if $M$ is a $K$-quasiconformally homogeneous planar domain which is not simply connected,  then $K\geq 1+\epsilon$, \cite{KM}.  Earlier work looked at this problem on hyperbolic manifolds in dimension $n\geq 3$  see \cite{BT0}  and also in Riemann surfaces where the problem is still open though partial results can be found,   see \cite{BT1,BT2}.

\medskip

One of the reasons this paper is widely cited is a key tool identified and developed for the first time; namely the quasihyperbolic metric -- defined below -- and estimates establishing how it distorts under quasiconformal mappings.

\medskip

\noindent{\bf Quasihyperbolic geometry}   Subsequently  Gehring wrote  with his former student Brad Osgood,  {\em Uniform domains and the quasihyperbolic metric}, which appeared in   J. Analyse Math.,  36 (1979), 50 -- 74.  Uniform domains were introduced by O. Martio and J. Sarvas around  1979 in connection with the injectivity properties of functions \cite{MS}. They found their way into P W  Jones' work characterising the domains $\Omega\subset\IR^n$ for which each BMO function $u$ in $\Omega$ has a BMO extension to ${\IR}^n$.   A proper subdomain $\Omega \subset{\IR}^n$, $n\geq 2$, is said to be uniform if there are $0<a,b<\infty$ such that any pair of points $x_1,x_2\in \Omega$ can be joined by a rectifiable curve $\gamma\subset \Omega$ for which $s(\gamma)\leq a|x_1-x_2|$ and $\min_{j=1,2}s(\gamma(x_j,x))\leq bd(x,\partial D)$ for all $x\in\gamma$, where $s(\gamma)$ denotes the Euclidean length of $\gamma$, $\gamma(x_j,x)$ is the part of $\gamma$ between $x_j$ and $x$, and $d(x,\partial D)$ is the Euclidean distance from $x$ to $\partial D$. 

Simply connected planar domains are uniform if and only if they are quasidisks. 

Again  the quasihyperbolic metric,  $k_D$  for a domain $
D$,  is key tool.  This metric is defined by 
\[ k_D(x_1,x_2)=\inf\int_\gamma d(x,\partial D)^{-1}\,ds \]
 where the infimum is taken over all rectifiable curves $\gamma\subset D$ joining $x_1$ to $x_2$ in $D$. It is not difficult  to see that 
 \[ k_D(x_1,x_2)\geq j_D(x_1,x_2)={\textstyle\frac 1{2}}\log\Big[\Big(\frac{|x_1-x_2|}{d(x_1,\partial D)}+1\Big)\Big(\frac{|x_1-x_2|}{d(x_2,\partial D)}+1\Big)\Big] \]
  for all $x_1,x_2\in D$. Gehring and Osgood first prove using the geodesics of $k_D$ (which turn out to be $C^{1,1}$ curves) that $D$ is uniform if and only if there are constants $0<c,d<\infty$ such that $k_D(x_1,x_2)\leq cj_D(x_1,x_2)+d$ for all $x_1,x_2\in D$ and this was actually what P. Jones used in \cite{PJ1,PJ2}. They also proved that $k_D$ and $j_D$ are quasi-invariant under quasiconformal mappings of $\Omega$ and this gives another proof of the invariance of the family of uniform domains under quasiconformal mappings of $\overline{{\IR}}^n$; an earlier result of Olli Martio \cite{OM}.   Gehring and Osgood then define a domain $D\subset{\IR}^2$ to be quasiconformally decomposable if there is $K\in[1,\infty)$ such that for any $x_1,x_2\in D$ there is a quasidisk $D_0\subset D$ which contains $x_1$ and $x_2$. This gives a surprising new characterisation of uniform domains : $D\subset{\IR}^2$ is uniform if and only if $D$ is quasiconformally decomposable, and via this characterization they give an alternative proof for the main injectivity properties of uniform domains in ${\IR}^2$ that were discovered by  Martio and J. Sarvas \cite{MS}.  This decomposition theorem is false in higher dimensions $n\geq 3$,  but is true if one redefines decomposability via the equivalent assumption in two dimensions,  that  there is $L\geq 1$ such that every pair of points lies in the $L$-bilipschitz image of a ball,  \cite{MartinThesis}.
 
\medskip

\noindent{\bf Discrete groups and hyperbolic geometry. } In the last couple of decades of Gehring's research career he and G J  Martin worked on the geometry of discrete groups,  writing about 30 joint papers.  This represented a new research direction for Gehring and was motivated by his attendance at Alan Beardon's  series of lectures on the geometry of discrete groups in Ann Arbor around 1980 - Beardon was just putting the finishing touches on his excellent book \cite{Beardon}.  Gehring  started thinking about how the general theory of M\"obius (or $1$-quasiconformal) groups would generalise to arbitrary discrete groups of quasiconformal mappings - typically with the assumption of a uniform bound on the distortion and touched on in his homogeneity paper with Palka.   It became clear that a surprising amount of the general theory of discrete groups used only the compactness properties of quasiconformal mappings.  These initial studies became the paper {\em Discrete quasiconformal groups. I.},  Proc. London Math. Soc.,  55 (1987),  331--358.  Part II never was finished.  The theory of convergence groups, identified and developed in that paper, was a timely idea and meshed well with Mikhael Gromov's theory of hyperbolic groups \cite{Gromov} published in the same year.  In the hands of Pekka Tukia it was quickly used   to characterise certain conjugates of M\"obius groups among groups of homeomorphisms \cite{Tukia}, a programme that was completed by David Gabai \cite{Gabai} and Andrew Casson \& Douglas Jungreis \cite{CJ} and,  with earlier work of Geoffery Mess,  led directly to the resolution of the Seifert conjecture :   a compact orientable irreducible $3$-manifold with  fundamental group having an infinite cyclic normal subgroup is Seifert fibered. Both results include a far reaching generalisation of the Nielsen Realisation problem,  earlier solved by Steve Kerckhoff \cite{SK}.  Mike Freedman even showed an equivalence  between the four-dimensional surgery conjecture and an extension problem for convergence groups of $\oR^3$,  \cite{Freedman}.   These new and surprising connections between quasiconformal mappings and low dimensional topology and so forth was enormously exciting to Gehring,  and some of the most interesting unresolved questions in low dimensional topology - for instance the Cannon Conjecture,  a group-theoretic generalization of the generic case of W. P. Thurston's famous Geometrization Conjecture  --  can be framed within the theory of convergence groups,  see \cite{Can} and \cite{MS}.

A group of homeomorphisms acting on $\oR^n$ is called a {\em quasiconformal group} if there is some finite $K$ such that each $g\in G$ is $K$-quasiconformal.  In two dimensions every quasiconformal group acting on $\oR^2=\oC$,  the Riemann sphere, is the quasiconformal conjugate of a conformal group or a M\"obius group as result of Dennis Sullivan and proved by  Tukia \cite{Tukia2},  that is there is a quasiconformal $f:\oC\to\oC$ and a group of M\"obius transformations $\Gamma$ such that $G=f\Gamma f^{-1}$. On the other hand, if $n\geq 3$, then there is a discrete quasiconformal group not quasiconformally conjugate to a M\"obius group \cite{Tukia3,MartinDQC}.  
 
 A {\em convergence group} $G$ is  a group of self homeomorphisms of $\oR^n$ which has the following property: every infinite subfamily of $G$ contains a sequence $\{g_j\}$ such that one of the following holds: 
 \begin{enumerate}
 \item there is a self-homeomorphism $g$ of $\oR^n$ such that $g_j\rightarrow g$ and $g_j^{-1}\rightarrow g^{-1}$ uniformly in $\oR^n$ as $j\rightarrow\infty$, 
 \item  there are points $x_0$ and $y_0$ in $\oR^n$ such that $g_j\rightarrow y_0$ and $g_j^{-1}\rightarrow x_0$ locally uniformly in $\oR^n-\{x_0\}$ and $\oR^n-\{y_0\}$, respectively, as $j\rightarrow\infty$. 
 \end{enumerate} 
 Convergence groups are similarly defined on other spaces.  Quasiconformal groups are convergence groups. A convergence group is {\em discrete} if (1) above never holds. The limit set $L(G)$ of a discrete convergence group $G$ is  defined  as for discrete M\"obius groups - the points of accumulation of a generic orbit. It has quite similar structure to that of a discrete M\"obius group:  $L(G)$ is either nowhere dense or is equal to $\oR^n$.  If $L(G)$ contains three points, then $L(G)$ is a perfect set. The elements of a discrete convergence group fall  into three kinds, elliptic, parabolic and loxodromic.  Examples are obtained from the following: if $E$ is a totally disconnected closed set in $\oR^n$ and if a group $G$ of self-homeomorphisms of $\oR^n$ is properly discontinuous in $\oR^n\setminus E$, then $G$ is a discrete convergence group whose limit set lies in $E$.
 
This initial foray into the geometry of discrete groups led to further research in {\em Kleinian groups} motivated by attempts to generalize to higher dimensions the universal geometric constraints on Fuchsian groups that  Beardon gave an exposition of \cite{Beardon}.  They observed a very interesting connection between criteria for discrete groups and holomorphic dynamics - the iteration of polynomials of one complex variable on the Riemann sphere,  \cite{GM1}.  This connection was already inherent in J\o rgensen's proof of his remarkable inequality for discrete groups \cite{Jor}. This polynomial came from a trace identity in $SL(2,\IC)$,  basically the Fricke identity.  Gehring and Martin had come to Mittag-Leffler in 1990 to work on this project as part of a six month thematic programme.  They discovered entirely new classes of polynomial trace identities and thus an entirely new tool to study the geometry of Kleinian groups which led to a number of new developments.  Among the more important papers is {\em Commutators, collars and the geometry of M\"obius groups},  J. Anal. Math.,  63 (1994), 175 -- 219.   Jane Gilman wrote  {\em 
``This paper is the seventh in a series of ten remarkable papers that Gehring and Martin have written on the geometry of discrete M\"obius groups . . . The papers use a combination of hyperbolic geometry and complex iteration theory to obtain discreteness criteria for M\"obius groups.  The results obtained in the current paper include a sharp analogue (for subgroups of $PSL(2,\IC)$ with an elliptic element) to the Shimizu-Leutbecher inequality  . . . and the elimination of a large region of possible values for the commutator parameter (trace$[g,h] - 2$ with $g,h \in G$) of a discrete group $G$ with two elliptic generators. From this the authors are able to obtain a sharp lower bound for the distance between the axes of elliptics of order $n$ and $m$ in any discrete groups (for many values of $m$ and $n$). The results have as a corollary substantial improvements ... on volumes of hyperbolic orbifolds  ... .''} Even with these new tools at hand it took quite some time to claim the real prize.  This came in the paper {\em Minimal co-volume hyperbolic lattices. I. The spherical points of a Kleinian group},  
Ann. of Math., {\bf 170} (2009),  123--161.   Gehring's   last research article.  This paper, along with a sequel \cite{MM}  completes the identification of the hyperbolic 3-orbifold of smallest volume. This solves a problem of Siegel from 1942  in three dimensions, also Problem 3.60 (F) in the Kirby problem list.   
  
  \medskip
 
\noindent{\bf Quasidisks.}  Finally we should talk about the ubiquitous quasidisk.  A {\em quasidisk} $\Omega$ is the image of the unit disk $\ID$ under a quasiconformal mapping $f:\IC\to\IC$,  so $f(\ID)=\Omega$.  Thus a quasidisk is a simply connected planar domain with reasonable geometric control on the structure of its boundary by way of such things as Hausdorff dimension, and ``turning conditions''.  The importance of the concept is reflected in the many remarkable and diverse applications in planar geometric function theory and also in low-dimensional geometry and topology where quasidisks form the components of ordinary set (where the action is properly discontinuous)  quasifuchsian groups.  In dynamics,  quasidisks form the components of the filled in Julia set of a hyperbolic rational map.  Further, the theory of holomorphic motions shows that when suitably interpreted a holomophic perturbation of the unit disk in the space of injections gives a quasidisk.  
  It was another lifelong task of Gehring to collect characterisations of quasidisks,  to give new and different proofs for these,  to find new applications of the theory and to generally spread the word about these wonderful objects.
 
 The nearest thing to a book that Gehring wrotehimself and in his lifetime is the monograph {\em Characteristic properties of quasidisks} which appears in  
S\'{e}minaire de Math\'{e}matiques Sup\'{e}rieures, 84. Presses de l'Universit\'{e} de Montr\'{e}al,  1982, 97 pp.   This book was subsequently expanded into the long awaited book {\em The ubiquitous quasidisk},  
in Mathematical Surveys and Monographs, 184. AMS,  2012  with Kari Hag. 
 \begin{center}  
\scalebox{0.32}{\includegraphics*{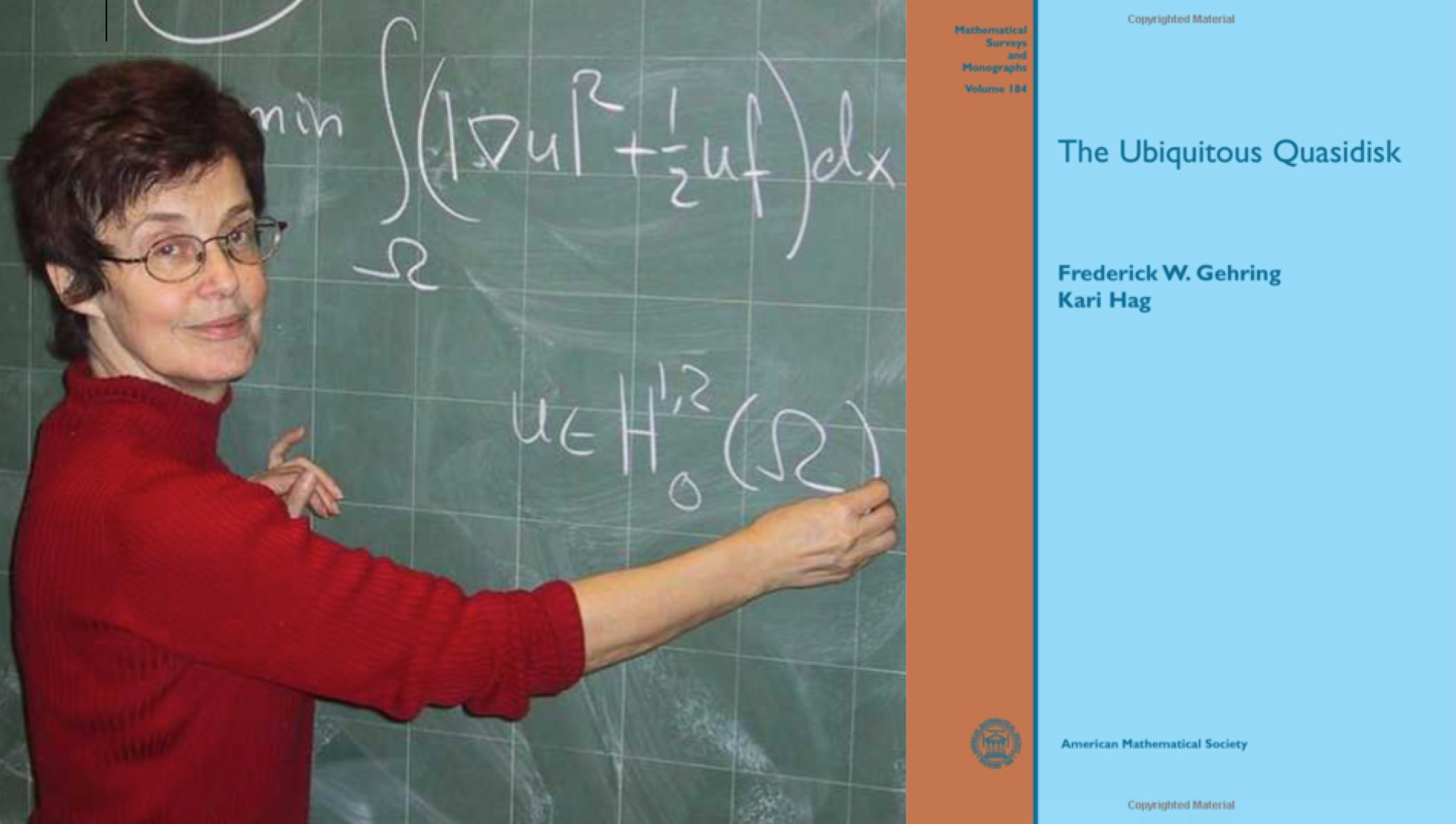}}\\
Former student Kari Hag and  book with Gehring
 \end{center}

These books can best thought of as a cross section of the many connections between the theory of planar quasiconformal mappings,  geometric function theory and analysis.  In the first set of notes,  Gehring gave 17 different characterisations,  this was subsequently extended to 30 (!).    Ahlfors' criterion  has been known for rather a long time; the Jordan domain $\Omega$ is a quasidisk if and only if there exists a constant $c$ such that, for all pairs of points 
\[ z_1,z_2\in\partial \Omega, \hskip10pt \min_{j=1,2}\; \big\{ \text{diam}(\gamma_j)\leq dc\, |z_1-z_2| \big\}\]
 where $\gamma_1$, $\gamma_2$ are the components of $\partial \Omega\smallsetminus\{z_1,z_2\}$. This is typical of the geometric charaterisations found,  others include the notions of linear local connectedness and so forth. Other geometric characterisations are based on the estimates on the hyperbolic or quasihyperbolic metric and similar sorts of things.  

Function-theoretic characterisation include the  Schwarzian derivative property.  For an analytic $\varphi$
\[ S\varphi = \Big(\frac{\varphi''}{\varphi'}\Big)'-\frac{1}{2}\Big(\frac{\varphi''}{\varphi'}\Big)^2, \]
and $\Omega$ is a quasidisk if and only if there exists a constant $c>0$ such that $f$ is injective whenever $f$ is analytic in $\Omega$ with $|S_\varphi|\leq c \rho_\Omega^2$, $\varphi'\neq 0$ in $\Omega$. Of course the Schwarzian derivative is  invariant under M\"obius transformations and has wide application in  modular forms and hypergeometric functions and Teichm\"uller theory.   Further criteria found by Gehring and Hag  concern  the injectivity properties of quasi-isometries, extension properties of  BMO or $W^{1,2}$ spaces (and their near relatives). 
A lot of the discussion concerns the many  characterizations involving the hyperbolic  and the quasihyperbolic metric.     We recommend the reader to thumb through the book to see the vast and diverse collection of results which cover a broad spectrum of planar function theory.

\newpage

\noindent{\bf Frederick W. Gehring: Published Works}.

\bigskip

\noindent Gehring, F. W., {\em Images of convergent sequences in sets}, J. London Math. Soc., {\bf  26}, (1951), 249--256.

 \noindent Burkill, J. C.; Gehring, F. W., {\em A scale of integrals from Lebesgue's to Denjoy's}, Quart. J. Math., Oxford Ser., {\bf 4}, (1953), 210--220. \\

 \noindent Gehring, F. W., {\em A study of $\alpha$-variation. I}, Trans. Amer. Math. Soc., {\bf  76}, (1954), 420--443. \\

 \noindent  Gehring, F. W., {\em A note on a paper by L. C. Young}, Pacific J. Math., {\bf  5}, (1955), 67--72.  \\

 \noindent  Gehring, F. W., {\em On the Dirichlet problem}, Michigan Math. J., {\bf  3}, (1955--1956), 201. \\

 \noindent  Gehring, F. W., {\em On the radial order of subharmonic functions}, J. Math. Soc. Japan, {\bf  9}, (1957), 77--79. \\

 \noindent  Gehring, F. W., {\em The Fatou theorem and its converse}, Trans. Amer. Math. Soc., {\bf  85}, (1957), 106--121.  \\

  \noindent  Gehring, F. W., {\em The Fatou theorem for functions harmonic in a half-space}, Proc. London Math. Soc., {\bf 8}, (1958), 149--160.  \\

  \noindent  Gehring, F. W., {\em The asymptotic values for analytic functions with bounded characteristic},  Quart. J. Math. Oxford Ser., {\bf 9}, (1958),  282--289.  \\

  \noindent  Gehring, F. W.; Lohwater, A. J., {\em On the Lindel\"of theorem},  Math. Nachr., {\bf  19}, (1958), 165--170.  \\

  \noindent  Gehring, F. W., {\em On solutions of the equation of heat conduction},  Michigan Math. J., {\bf  5}, (1958), 191--202.  \\

 \noindent   Gehring, F. W.; Lehto, O., {\em On the total differentiability of functions of a complex variable}, Ann. Acad. Sci. Fenn. Ser. A I, {\bf 272}, (1959),  9 pp.  \\

 \noindent   Gehring, F. W., {\em  The boundary behavior and uniqueness of solutions of the heat equation},  Trans. Amer. Math. Soc., {\bf  94}, (1960), 337--364.  \\

 \noindent    Gehring, F. W., {\em  Harmonic functions and Tauberian theorems},  Proc. London Math. Soc., {\bf 10}, (1960), 88--106. \\

  \noindent  Gehring, F. W., {\em  The definitions and exceptional sets for quasiconformal mappings},  Ann. Acad. Sci. Fenn. Ser. A I, {\bf 281}, (1960), 28 pp. \\

  \noindent  Gehring, F. W.; Haahti, H., {\em  The transformations which preserve the harmonic functions},  Ann. Acad. Sci. Fenn. Ser. A I, {\bf 293}, (1960), 12 pp.  \\

  \noindent  Gehring, F. W., {\em  Rings and quasiconformal mappings in space},  Proc. Nat. Acad. Sci. U.S.A., {\bf  47}, (1961), 98--105.  \\

  \noindent  Gehring, F. W., {\em  Symmetrization of rings in space}, Trans. Amer. Math. Soc., {\bf  101}, (1961), 499--519.  \\

  \noindent  Gehring, F. W., {\em  A remark on the moduli of rings},  Comment. Math. Helv., {\bf  36}, (1961), 42--46.  \\

  \noindent  Gehring, F. W.; V\"ais\"al\"a, J., {\em On the geometric definition for quasiconformal mappings},  Comment. Math. Helv., {\bf  36}, (1961), 19--32.  \\

  \noindent  Gehring, F. W. {\em Rings and quasiconformal mappings in space},  Trans. Amer. Math. Soc., {\bf  103}, (1962), 353--393.  \\

  \noindent   Gehring, F. W., {\em  Extremal length definitions for the conformal capacity of rings in space},  Michigan Math. J., {\bf  9}, (1962), 137--150. \\

  \noindent   Gehring, F. W.; Hayman, W. K., {\em  An inequality in the theory of conformal mapping},  J. Math. Pures Appl., {\bf 41}, (1962), 353--361. \\

   \noindent   Gehring, F. W., {\em  Quasiconformal mappings in space},  Bull. Amer. Math. Soc., {\bf  69},  (1963), 146--164. \\

 \noindent    Gehring, F. W.; af H\"allstr\"om, G., {\em  A distortion theorem for functions univalent in an annulus},  Ann. Acad. Sci. Fenn. Ser. A I , {\bf  325},  (1963), 16 pp.\\

  \noindent   Gehring, F. W., {\em  The Carath\'eodory convergence theorem for quasiconformal mappings in space},  Ann. Acad. Sci. Fenn. Ser. A I , {\bf 336}, (1963), 21 pp. \\

  \noindent   Agard, S. B.; Gehring, F. W., {\em  Angles and quasiconformal mappings},  Proc. London Math. Soc., {\bf 14}, (1965), 1--21. \\

  \noindent   Gehring, F. W., {\em  Extension of quasiconformal mappings in three space},  J. Analyse Math., {\bf  14 }, (1965), 171--182. \\

   \noindent  Gehring, F. W.; V\"ais\"al\"a, J., {\em  The coefficients of quasiconformality of domains in space},  Acta Math., {\bf  114}, (1965), 1--70. \\

   \noindent  Gehring, F. W.; Reich, E., {\em  Area distortion under quasiconformal mappings},  Ann. Acad. Sci. Fenn. Ser. A I, {\bf 388}, (1966), 15 pp. \\

 \noindent    Gehring, F. W., {\em  Coefficients of quasiconformality of domains in three space},  1966 Contemporary Problems in Theory Anal. Functions (Internat. Conf., Erevan, 1965)   pp. 83--88 Izdat. ``Nauka'', Moscow. \\

 \noindent    Gehring, F. W., {\em  Definitions for a class of plane quasiconformal mappings},  Nagoya Math. J., {\bf  29}, (1967), 175--184. \\

 \noindent    Gehring, F. W., {\em  Extension theorems for quasiconformal mappings in $n$-space},  J. Analyse Math., {\bf  19}, (1967), 149--169.\\

 \noindent    Gehring, F. W., {\em  Extension theorems for quasiconformal mappings in $n$-space},  1968 Proc. Internat. Congr. Math. (Moscow, 1966) pp. 313--318 Izdat. ``Mir'', Moscow.\\

  \noindent   Gehring, F. W., {\em  Quasiconformal mappings of slit domains in three space},  J. Math. Mech., {\bf  18}, (1969), 689--703.\\

 \noindent    Gehring, F. W., {\em  Quasiconformal mappings which hold the real axis pointwise fixed},  1970, Mathematical Essays Dedicated to A. J. Macintyre, pp. 145--148, Ohio Univ. Press, Athens, Ohio. \\

  \noindent  Gehring, F. W., {\em  Extremal mappings of tori},   Certain problems of mathematics and mechanics (on the occasion of the seventieth birthday of M. A. Lavrentiev)  pp. 146--152. Izdat. "Nauka'', Leningrad, 1970. \\

 \noindent   Gehring, F. W., {\em  Inequalities for condensers, hyperbolic capacity, and extremal lengths},  Michigan Math. J., {\bf  18}, (1971), 1--20.\\

 \noindent   Gehring, F. W., {\em  Lipschitz mappings and the $p$-capacity of rings in $n$-space}, Advances in the theory of Riemann surfaces (Proc. Conf., Stony Brook, N.Y., 1969), pp.175--193. Ann. of Math. Studies, {\bf 66}. Princeton Univ. Press, Princeton, N.J., 1971. \\

  \noindent    Gehring, F. W.; Huckemann, F., {\em  Quasiconformal mappings of a cylinder},  Proc. Romanian-Finnish Seminar on Teichm\"uller Spaces and Quasiconformal Mappings (Bra\c{s}ov, 1969), pp.171--186. (Publ. House of the Acad. of the Socialist Republic of Romania, Bucharest, 1971)\\

  \noindent   Gehring, F. W., {\em  Dilatations of quasiconformal boundary correspondences},  Duke Math. J., {\bf  39}, (1972), 89--95. \\

  \noindent  Gehring, F. W., {\em  The $L^p$-integrability of the partial derivatives of quasiconformal mapping},  Bull. Amer. Math. Soc., {\bf  79}, (1973), 465--466. \\

 \noindent   Gehring, F. W.; V\"ais\"al\"a, J., {\em  Hausdorff dimension and quasiconformal mappings},  J. London Math. Soc., {\bf 6}, (1973), 504--512.  \\

  \noindent  Gehring, F. W., {\em  The $L^p$-integrability of the partial derivatives of a quasiconformal mapping},  Acta Math., {\bf  130}, (1973), 265--277.  \\

 \noindent   Gehring, F. W.; Kelly, J. C., {\em  Quasiconformal mappings and Lebesgue density},  Discontinuous groups and Riemann surfaces (Proc. Conf., Univ. Maryland, College Park, Md., 1973), pp. 171--179. Ann. of Math. Studies,  {\bf  79}, Princeton Univ. Press, Princeton, N.J., 1974. \\

 \noindent   Gehring, F. W., {\em  The Hausdorff measure of sets which link in Euclidean space},  Contributions to analysis (a collection of papers dedicated to Lipman Bers), pp. 159--167. Academic Press, New York, 1974.  \\

  \noindent  Gehring, F. W., {\em  Inequalities for condensers, hyperbolic capacity, and extremal lengths},  Topics in analysis (Colloq. Math. Anal., Jyv\"askyl\"a, 1970), pp.133--136. Lecture Notes in Math., Vol. 419, Springer, Berlin, 1974.  \\

 \noindent   Gehring, F. W., {\em  The $L^p$-integrability of the partial derivatives of a quasiconformal mapping},  Proceedings of the Symposium on Complex Analysis (Univ. Kent, Canterbury, 1973), pp. 73--74. London Math. Soc. Lecture Note Ser.,   {\bf 12}, Cambridge Univ. Press, 1974.  \\

 \noindent   Gehring, F. W., {\em  Lower dimensional absolute continuity properties of quasiconformal mappings},  Math. Proc. Cambridge Philos. Soc., {\bf  78}, (1975), 81--93.  \\

  \noindent  Gehring, F. W., {\em  Quasiconformal mappings in ${\mathbb R}^n$},  Lectures on quasiconformal mappings,  44--110. Dept. Math., Univ. Maryland, Lecture Note, {\bf 14}, Dept. Math., Univ. Maryland, College Park, Md., 1975.  \\

 \noindent   Gehring, F. W., {\em  Some metric properties of quasiconformal mappings},  Proceedings of the International Congress of Mathematicians (Vancouver, B. C., 1974), Vol. 2, pp. 203--206. Canad. Math. Congress, Montreal, Que., 1975.  \\

  \noindent  Gehring, F. W.; Palka, B. P., {\em  Quasiconformally homogeneous domains}, J. Analyse Math., {\bf  30}, (1976), 172--199.  \\

 \noindent   Gehring, F. W., {\em  Absolute continuity properties of quasiconformal mappings},  Symposia Mathematica, Vol. LXVIII (Convegno sulle Transformazioni Quasiconformi e Questioni Connesse, INDAM, Rome, 1974), pp. 551--559. Academic Press, London, 1976.   \\

  \noindent  Gehring, F. W., {\em  Quasiconformal mappings},  Complex analysis and its applications (Lectures, Internat. Sem., Trieste, 1975), Vol. II, pp. 213--268. Internat. Atomic Energy Agency, Vienna, 1976.    \\

 \noindent   Gehring, F. W., {\em  A remark on domains quasiconformally equivalent to a ball},  Ann. Acad. Sci. Fenn. Ser. A I Math., {\bf  2}, (1976), 147--155.  \\

 \noindent   Gehring, F. W., {\em  Univalent functions and the Schwarzian derivative},  Comment. Math. Helv., {\bf  52}, (1977), 561--572.  \\

  \noindent  Gehring, F. W., {\em  Spirals and the universal Teichm\"uller space},  Acta Math., {\bf   141}, (1978),  99--113.  \\

  \noindent  Gehring, F. W., {\em  Some problems in complex analysis},  Proceedings of the First Finnish-Polish Summer School in Complex Analysis (Podlesice, 1977), Part II, pp. 61--64, Univ. {\L}\'od\'z,  1978.  \\

  \noindent  Gehring, F. W., {\em  Remarks on the universal Teichm\"uller space},  Enseign. Math., {\bf 24}, (1978),   173--178.  \\

 \noindent   Gehring, F. W., {\em  Univalent functions and the Schwarzian derivative},  Proceedings of the Rolf Nevanlinna Symposium on Complex Analysis (Math. Res. Inst., Univ. Istanbul, Silivri, 1976), pp. 19--24, Publ. Math. Res. Inst. Istanbul, 7, Univ. Istanbul, Istanbul, 1978. \\

  \noindent  Gehring, F. W.; Osgood, B. G., {\em  Uniform domains and the quasihyperbolic metric},  J. Analyse Math., {\bf  36}, (1979), 50--74 (1980).  \\

  \noindent   Beardon, A. F.; Gehring, F. W., {\em  Schwarzian derivatives, the Poincar\'e metric and the kernel function},  Comment. Math. Helv., {\bf  55}, (1980), 50--64. \\

   \noindent  Gehring, F. W., {\em  Spirals and the universal Teichm\"uller space},  Riemann surfaces and related topics: Proceedings of the 1978 Stony Brook Conference (State Univ. New York, Stony Brook, N.Y., 1978), pp. 145--148, Ann. of Math. Stud., {\bf  97}, Princeton Univ. Press, Princeton, N.J., 1981. \\

  \noindent  Gehring, F. W.; Hayman, W. K.; Hinkkanen, A., {\em Analytic functions satisfying H\"older conditions on the boundary},  J. Approx. Theory, {\bf  35}, (1982), 243--249. \\

  \noindent  Teichm\"uller, Oswald Gesammelte Abhandlungen.  [Collected papers] Edited and with a preface by Lars V. Ahlfors and F W. Gehring. Springer-Verlag, Berlin-New York, 1982.   \\

 \noindent   Gehring, F. W., {\em  Characteristic properties of quasidisks},  S\'eminaire de Math\'ematiques Sup\'erieures [Seminar on Higher Mathematics], 84. Presses de l'Universit\'e de Montr\'eal, Montreal, Que., 1982. 97 pp.   \\

  \noindent  Gehring, F. W.  {\em Injectivity of local quasi-isometries},    Comment. Math. Helv., {\bf  57}, (1982),  202--220. \\

  \noindent  Garnett, J. B.; Gehring, F. W.; Jones, P. W., {\em  Conformally invariant length sums}, Indiana Univ. Math. J., {\bf  32}, (1983),  809--829. \\

  \noindent  Gehring, F. W.; Martio, O., {\em  Quasidisks and the Hardy-Littlewood property}, Complex Variables Theory Appl., {\bf  2}, (1983),  67--78. \\

 \noindent   Astala, K.; Gehring, F. W. {\em Injectivity criteria and the quasidisk},   Complex Variables Theory Appl., {\bf  3}, (1984),  45--54. \\

 \noindent   Gehring, F. W.; Pommerenke, Ch., {\em  On the Nehari univalence criterion and quasicircles}, Comment. Math. Helv., {\bf  59}, (1984),  226--242. \\

 \noindent   Astala, K.; Gehring, F. W., {\em  Quasiconformal analogues of theorems of Koebe and Hardy-Littlewood}, Michigan Math. J., {\bf  32}, (1985), 99--107.\\

 \noindent   Anderson, J. M.; Gehring, F. W.; Hinkkanen, A., {\em  Polynomial approximation in quasidisks}, Differential geometry and complex analysis, 75--86, Springer, Berlin, 1985. \\

  \noindent  Gehring, F. W.; Martio, O., {\em  Lipschitz classes and quasiconformal mappings}, Ann. Acad. Sci. Fenn. Ser. A I Math., {\bf  10}, (1985), 203--219. \\

 \noindent   Gehring, F. W.; Martio, O., {\em  Quasiextremal distance domains and extension of quasiconformal mappings}, J. Analyse Math., {\bf  45}, (1985), 181--206. \\

 \noindent   Gehring, F. W., {\em  Extension of quasi-isometric embeddings of Jordan curves}, Complex Variables Theory Appl., {\bf  5}, (1986), 245--263. \\

  \noindent   Astala, K.; Gehring, F. W., {\em  Injectivity, the BMO norm and the universal Teichm\"uller space},  J. Analyse Math., {\bf  46}, (1986), 16--57. \\

  \noindent  Gehring, F. W., {\em  Uniform domains and the ubiquitous quasidisk}, Jahresber. Deutsch. Math.-Verein., {\bf  89}, (1987),  88--103. \\

 \noindent   Gehring, F. W.; Martin, G. J., {\em  Discrete quasiconformal groups. I}, Proc. London Math. Soc., {\bf 55}, (1987),  331--358. \\

 \noindent   Gehring, F. W.; Martin, G. J., {\em  Discrete convergence groups}, Complex analysis, I (College Park, Md., 1985--86), 158--167, Lecture Notes in Math., 1275, Springer, Berlin, 1987. \\

 \noindent   Gehring, F. W.; Hag, K., {\em  Remarks on uniform and quasiconformal extension domains}, Complex Variables Theory Appl., {\bf  9}, (1987), 175--188. \\

  \noindent  Gehring, F. W., {\em  Topics in quasiconformal mappings}, Proceedings of the International Congress of Mathematicians, Vol. 1, 2 (Berkeley, Calif., 1986), 62--80, Amer. Math. Soc., Providence, RI, 1987.\\

  \noindent  Holomorphic functions and moduli. I \& II. Proceedings of the workshop held in Berkeley, California, March 13--19, 1986. Edited by D. Drasin, C. J. Earle, F. W. Gehring, I. Kra and A. Marden. Mathematical Sciences Research Institute Publications, 11. Springer-Verlag, New York, 1988.\\

 \noindent   Garnett, J. B.; Gehring, F. W.; Jones, P. W., {\em  Quasiconformal groups and the conical limit set}, Holomorphic functions and moduli, Vol. II (Berkeley, CA, 1986), 59--67, Math. Sci. Res. Inst. Publ., {\bf  11}, Springer, New York, 1988. \\

  \noindent  Gehring, F. W.; Martin, G. J., {\em  The matrix and chordal norms of M\"obius transformations}, Complex analysis, 51--59, Birkh\"auser, Basel, 1988. \\

  \noindent  Gehring, F W., {\em  Quasiconformal mappings}, A plenary address presented at the International Congress of Mathematicians held in Berkeley, California, August 1986. Introduced by G. D. Mostow. ICM Series. American Mathematical Society, Providence, RI, 1988. 1 videocassette (NTSC; 1/2 inch; VHS) (60 min.)\\

 \noindent   Gehring, F. W.; Martin, G. J., {\em  Iteration theory and inequalities for Kleinian groups}, Bull. Amer. Math. Soc. (N.S.), {\bf  21}, (1989),   57--63. \\

  \noindent  Gehring, F. W.; Pommerenke, Ch., {\em  Circular distortion of curves and quasicircles},  Ann. Acad. Sci. Fenn. Ser. A I Math., {\bf  14}, (1989),   381--390.\\

  \noindent  Gehring, F. W.; Martin, G. J., {\em  Stability and extremality in J\o rgensen's inequality},   Complex Variables Theory Appl., {\bf  12}, (1989),   277--282. \\

  \noindent  Gehring, F. W.; Hag, K.; Martio, O., {\em  Quasihyperbolic geodesics in John domains},   Math. Scand., {\bf  65}, (1989),  75--92. \\

  \noindent  Astala, K.; Gehring, F. W., {\em  Crickets, zippers and the Bers universal Teichm\"uller space},   Proc. Amer. Math. Soc., {\bf  110}, (1990),  675--687. \\

 \noindent   Halmos, Paul R.; Gehring, F. W., {\em  Allen L. Shields},   Math. Intelligencer, {\bf  12}, (1990),  20.\\

  \noindent  Gehring, F. W.; Martin, G. J., {\em  Inequalities for M\"obius transformations and discrete groups},   J. Reine Angew. Math., {\bf  418}, (1991), 31--76. \\

  \noindent  Paul Halmos., {\em  Celebrating 50 years of mathematics},   Edited by John H. Ewing and F. W. Gehring. Springer-Verlag, New York, 1991.   \\

  \noindent  Gehring, F. W.; Martin, G. J., {\em  Some universal constraints for discrete M\"obius groups},  Paul Halmos, 205--220, Springer, New York, 1991. \\

  \noindent  Gehring, F. W.; Martin, G. J., {\em  Axial distances in discrete M\"obius groups},   Proc. Nat. Acad. Sci. U.S.A., {\bf  89}, (1992),   1999--2001. \\

 \noindent   Gehring, F. W., {\em  Topics in quasiconformal mappings: Quasiconformal space mappings},   20--38, Lecture Notes in Math., 1508, Springer, Berlin, 1992. \\

  \noindent  Gehring, F. W.; Martin, G. J., {\em  6-torsion and hyperbolic volume},   Proc. Amer. Math. Soc., {\bf  117}, (1993), 727--735. \\

  \noindent  Gehring, F. W.; Martin, G. J., {\em  The iterated commutator in a Fuchsian group},   Complex Variables Theory Appl., {\bf  21}, (1993), 207--218. \\

  \noindent  Gehring, F. W.; Martin, G. J., {\em  On the minimal volume hyperbolic 3-orbifold},   Math. Res. Lett., {\bf  1}, (1994),  107--114. \\

  \noindent  Gehring, F. W.; Martin, G. J., {\em  Commutators, collars and the geometry of M\"obius groups},   J. Anal. Math., {\bf  63}, (1994), 175--219. \\

  \noindent  Gehring, F. W.; Martin, G. J., {\em  Chebyshev polynomials and discrete groups},   Proceedings of the Conference on Complex Analysis (Tianjin, 1992), 114--125, Conf. Proc. Lecture Notes Anal., I, Int. Press, Cambridge, MA, 1994. \\

  \noindent  Gehring, F. W.; Martin, G. J., {\em  Holomorphic motions, Schottky's theorem and an inequality for discrete groups},   Computational methods and function theory 1994 (Penang), 173--181, Ser. Approx. Decompos., {\bf  5}, World Sci. Publ., River Edge, NJ, 1995. \\

  \noindent  Gehring, F. W.; Martin, G. J., {\em  On the Margulis constant for Kleinian groups. I},  Ann. Acad. Sci. Fenn. Math., {\bf  21}, (1996),   439--462. \\

  \noindent  Gehring, F. W.; Maclachlan, C.; Martin, G. J., {\em  On the discreteness of the free product of finite cyclic groups},  Mitt. Math. Sem. Giessen, {\bf 228}, (1996), 9--15. \\

  \noindent  Gehring, F. W.; Maclachlan, C.; Martin, G. J.; Reid, A. W., {\em  Arithmeticity, discreteness and volume},   Trans. Amer. Math. Soc., {\bf  349}, (1997), 3611--3643.\\

  \noindent  Gehring, F. W.; Martin, G. J., {\em  Geodesics in hyperbolic $3$-folds},   Michigan Math. J., {\bf  44}, (1997),  331--343. \\

 \noindent   Cao, C.; Gehring, F. W.; Martin, G. J., {\em  Lattice constants and a lemma of Zagier},   Lipa's legacy (New York, 1995), 107--120, Contemp. Math., {\bf  211}, Amer. Math. Soc., Providence, RI, 1997. \\

  \noindent  Gehring, F. W.; Martin, G. J., {\em  Hyperbolic axes and elliptic fixed points in a Fuchsian group},   Complex Variables Theory Appl., {\bf  33}, (1997), 303--309. \\

 \noindent   Gehring, F; Kra, I.; Osserman, R., {\em The mathematics of Lars Valerian Ahlfors},   Edited by Steven G. Krantz. Notices Amer. Math. Soc., {\bf  45}, (1998),  233--242. \\

 \noindent   Gehring, F. W.; Maclachlan, C.; Martin, G. J., {\em  Two-generator arithmetic Kleinian groups. II},   Bull. London Math. Soc., {\bf  30}, (1998),  258--266. \\

  \noindent  Gehring, F. W.; Marshall, T. H.; Martin, G. J., {\em  The spectrum of elliptic axial distances in Kleinian groups},   Indiana Univ. Math. J., {\bf  47}, (1998),   1--10. \\

 \noindent   Gehring, F. W.; Martin, G. J., {\em  Precisely invariant collars and the volume of hyperbolic $3$-folds},   J. Differential Geom., {\bf  49}, (1998),   411--435. \\

  \noindent  Gehring, F. W.; Marshall, T. H.; Martin, G. J., {\em  Collaring theorems and the volumes of hyperbolic $n$-manifolds},   New Zealand J. Math., {\bf 27}, (1998),   207--225. \\

  \noindent  Gehring, F. W.; Martin, G. J., {\em  The volume of hyperbolic $3$-folds with p-torsion, $p\geq 6$},   Quart. J. Math. Oxford Ser., {\bf  50}, (1999),  1--12. \\

 \noindent   Gehring, F. W.; Iwaniec, T., {\em  The limit of mappings with finite distortion},   Ann. Acad. Sci. Fenn. Math., {\bf  24}, (1999),   253--264. \\

 \noindent   Gehring, F. W.; Hag, K., {\em  Hyperbolic geometry and disks},   Continued fractions and geometric function theory (CONFUN) (Trondheim, 1997). J. Comput. Appl. Math., {\bf  105}, (1999),  275--284. \\

 \noindent   Gehring, F W.; Hag, K., {\em A bound for hyperbolic distance in a quasidisk},   Computational methods and function theory 1997 (Nicosia), 233--240, Ser. Approx. Decompos., {\bf  11}, World Sci. Publ., River Edge, NJ, 1999. \\

 \noindent   Gehring, F W., {\em  Characterizations of quasidisks},   Quasiconformal geometry and dynamics (Lublin, 1996), 11--41, Banach Center Publ., {\bf 48}, Polish Acad. Sci., Warsaw, 1999. \\

  \noindent  Gehring, F W., {\em  Variations on a theorem of Fej\'er and Riesz},   XII-th Conference on Analytic Functions (Lublin, 1998). Ann. Univ. Mariae Curie-Sk{\l}odowska Sect. A, {\bf 53}, (1999), 57--66. \\

 \noindent   Gehring, F W.; Hag, K., {\em Bounds for the hyperbolic distance in a quasidisk},   XII-th Conference on Analytic Functions (Lublin, 1998). Ann. Univ. Mariae Curie-Sk{\l}odowska Sect. A {\bf 53},(1999), 67--72.\\

  \noindent  Gehring, F. W.; Hag, K., {\em  The Apollonian metric and quasiconformal mappings},   In the tradition of Ahlfors and Bers (Stony Brook, NY, 1998), 143--163, Contemp. Math.,   {\bf 256},  Amer. Math. Soc., Providence, RI, 2000. \\

  \noindent  Gehring, F. W., {\em  The Apollonian metric},   Travaux de la Conf\'erence Internationale d'Analyse Complexe et du 8e S\'eminaire Roumano-Finlandais (Iassy, 1999). Math. Rep. (Bucur.) , {\bf 2},  (2000),  461--466 (2001).\\

  \noindent  Gehring, F. W.; Gilman, J. P.; Martin, G. J., {\em  Kleinian groups with real parameters},   Commun. Contemp. Math., {\bf  3}, (2001),   163--186. \\

  \noindent   Brown, E. H.; Cohen, F. R.; Gehring, F. W.; Miller, H. R.; Taylor, B. A., {\em  Franklin P. Peterson (1930--2000)},   Notices Amer. Math. Soc., {\bf  48}, (2001),  1161--1168. \\

 \noindent  Gehring, F. W.; Hag, K., {\em  Reflections on reflections in quasidisks},   Papers on analysis, 81--90, Rep. Univ. Jyv\"askyl\"a Dep. Math. Stat., {\bf 83},  Jyv\"askyl\"a, 2001. \\

  \noindent  Gehring, F W.; Marshall, Timothy H.; Martin, G J., {\em  Recent results in the geometry of Kleinian groups},   Comput. Methods Funct. Theory {\bf 2}, (2002),   249--256. \\

  \noindent  Gehring, F W.; Hag, K., {\em Sewing homeomorphisms and quasidisks},   Comput. Methods Funct. Theory {\bf 3}, (2003),  143--150. \\

 \noindent   Gehring, F. W., {\em  Quasiconformal mappings in Euclidean spaces},   Handbook of complex analysis: geometric function theory.  {\bf 2}, 1--29, Elsevier, Amsterdam, 2005. \\

  \noindent  Gehring, F. W.; Martin, G. J., {\em (p,q,r)-Kleinian groups and the Margulis constant},   Complex analysis and dynamical systems II, 149--169, Contemp. Math., {\bf 382}, Amer. Math. Soc., Providence, RI, 2005. \\

  \noindent  Gehring, F W.; Martin, G J., {\em Minimal co-volume hyperbolic lattices. I. The spherical points of a Kleinian group},   Ann. of Math. (2) {\bf 170}, (2009),  123--161. \\

 \noindent  Gehring, F W.; Hag, K., The ubiquitous quasidisk. With contributions by Ole Jacob Broch. Mathematical Surveys and Monographs, {\bf 184}.  American Mathematical Society, Providence, RI, 2012.  \\

\noindent G.J. Martin. Massey University,   New Zealand  
 

\begin{thebibliography}{99}  
   
 \bibitem{A}  L. V. Ahlfors, {\em On quasiconformal mappings},  J. Analyse Math., {\bf 3}, (1954), 1--58.
 
 \bibitem{AB2}  L. V. Ahlfors and A.  Beurling, {\em Invariants conformes et probl\`emes extr\'emaux}, C. R. Dixime Congrs Math. Scandinaves 1946, pp. 341--351. Jul. Gjellerups Forlag, Copenhagen, 1947.
 
 \bibitem{And} G.D.  Anderson, {\em The coefficients of quasiconformality of ellipsoids},  Ann. Acad. Sci. Fenn. Ser. A I No. 411,  (1967) 14 pp.
 
 \bibitem{Astala} K. Astala, {\em Area distortion of quasiconformal mappings},  Acta Math., {\bf 173},  (1994),   37--60.
 \bibitem{AIM} K. Astala, T.  Iwaniec and G.J. Martin, {\em Elliptic partial differential equations and quasiconformal mappings in the plane},  Princeton Mathematical Series, 48. Princeton University Press, Princeton, NJ, 2009.
 
 \bibitem{Beardon} A. F. Beardon, {\em The geometry of discrete groups},  Graduate Texts in Mathematics, 91. Springer-Verlag, New York, 1983.t
 
 \bibitem{Bo} B. V. Bojarski, {\em Generalized solutions of a system of differential equations of first order and of elliptic type with discontinuous coefficients}, 
Mat. Sb. N.S., {\bf  43}, (1957), 451--503.

\bibitem{BT0} P. Bonfert-Taylor, R. D. Canary,  G.J. Martin and E. Taylor, {\em Quasiconformal homogeneity of hyperbolic manifolds},  Math. Ann., {\bf  331},  (2005),  281--295.
\bibitem{BT1} P. Bonfert-Taylor,  G.J. Martin, A.W. Reid and E. Taylor, {\em Teichm\"uller mappings, quasiconformal homogeneity, and non-amenable covers of Riemann surfaces},  Pure Appl. Math. Q. 7 (2011),  Special Issue: In honor of Frederick W. Gehring, Part 2, 455 --468.
\bibitem{BT2} P. Bonfert-Taylor, Petra; R.D. Canary, J. Souto, and E.C. Taylor, {\em Exotic quasi-conformally homogeneous surfaces}, Bull. Lond. Math. Soc., {\bf 43},  (2011), 57--62.

\bibitem{Call} E.D. Callender, {\em 
H\"older continuity of n-dimensional quasi-conformal mappings}, 
Pacific J. Math., {\bf 10},  1960, 499--515. 

\bibitem{Can} J.W. Cannon, {\em The combinatorial Riemann mapping theorem}, Acta Math., {\bf 173}, (1994), 155-- 234 
 
 \bibitem{CJ} A. Casson and D. Jungreis, {\em 
Convergence groups and Seifert fibered 3-manifolds},  
Invent. Math., {\bf 118},  (1994), 441--456. 

\bibitem{Freedman} M.H. Freedman, {\em A geometric reformulation of 4-dimensional surgery}, Special volume in honor of R. H. Bing (1914--1986). Topology Appl., {\bf 24} (1986),  133--141.

\bibitem{GM1} F.W. Gehring and G.J. Martin, {\em  Iteration theory and inequalities for Kleinian groups},  Bull. Amer. Math. Soc., {\bf  21},  (1989),   57--63.

\bibitem{GM2} F.W. Gehring and G.J. Martin, {\em  Precisely invariant collars and the volume of hyperbolic 3-folds},  J. Differential Geom., {\em  49},  (1998),  411--435.

\bibitem{GM3}  F.W. Gehring, C. Maclachlan, G.J. Martin and A.W. Reid, {\em Arithmeticity, discreteness and volume}, Trans. Amer. Math. Soc., {\bf 349}, (1997),  3611--3643.

\bibitem{DG} E. De Giorgi, {\em Sulla differenziabilitˆ e l'analiticitˆ delle estremali degli integrali multipli regolari},  Mem. Accad. Sci. Torino. Cl. Sci. Fis. Mat. Nat., {\bf  3},  1957 25--43.

\bibitem{Gabai} D. Gabai, {\em Convergence groups are Fuchsian groups},  Ann. of Math., {\bf 136}, (1992),   447--510. 


\bibitem{Gromov} M. Gromov, {\em Hyperbolic groups},  Essays in group theory, 75--263, Math. Sci. Res. Inst. Publ., 8, Springer, New York, 1987.
 

\bibitem{Hag1}  K. Hag  and M. N\"a\"at\"anen, {\em On the outer coefficient of quasiconformality of a cylindrical map of a convex dihedral wedge},  Ann. Acad. Sci. Fenn. Ser. A I Math., {\bf  5}, (1980),  125 -- 130. 

\bibitem{Hag2}  K. Hag and M.K. Vamanamurthy, {\em The coefficients of quasiconformality of cones in n-space},  Ann. Acad. Sci. Fenn. Ser. A I Math., {\bf  3},  (1977),  267--275.  

\bibitem{HK} J. Heinonen and P. Koskela, {\em Quasiconformal maps in metric spaces with controlled geometry},  Acta Math., {\bf  181},  (1998),   1--61.

\bibitem{IwaniecGL} T. Iwaniec,  {\em The Gehring lemma},  Quasiconformal mappings and analysis (Ann Arbor, MI, 1995), 181Ð204, Springer, New York, 1998. 

 \bibitem{I1} T. Iwaniec, {\em $p$-harmonic tensors and quasiregular mappings}, Ann. of Math., {\bf 136},  (1992),  589--624.
 
 \bibitem{I2}  T. Iwaniec and G. Martin,  {\em Quasiregular mappings in even dimensions},  Acta Math., {\bf 170}, (1993),  29--81.
 \bibitem{IM}  T. Iwaniec and G. Martin,   {\em Geometric function theory and non-linear analysis},  Oxford Mathematical Monographs. The Clarendon Press, Oxford University Press, New York, 2001. 
 
 \bibitem{Jor} T. J\o rgensen, {\em On discrete groups of M\"obius transformations},  Amer. J. Math., {\bf 98},  (1976),  739--749.
 
 \bibitem{lehto}   O. Lehto, {\em Gehring Gehring and Finnish mathematics,  } Quasiconformal mappings and analysis. A collection of papers honoring F. W. Gehring. Papers from the International Symposium held in Ann Arbor, MI, August 1995. Edited by Peter Duren, Juha Heinonen, Brad Osgood and Bruce Palka. Springer-Verlag, New York, 1998. 

\bibitem{Lewis} L.G. Lewis, {\em Quasiconformal mappings and Royden algebras in space}, Trans. Amer. Math. Soc., {\bf 158},  1971 481--492.  

\bibitem{MM}  T.H. Marshall and G.J. Martin, {\em Minimal co-volume hyperbolic lattices, II: Simple torsion in a Kleinian group},  Ann. of Math., {\bf  176},  (2012),  261--301.

\bibitem{MartinThesis} G.J. Martin, {\em Quasiconformal and bi-Lipschitz homeomorphisms, uniform domains and the quasihyperbolic metric},  Trans. Amer. Math. Soc., {\em  292}, (1985),  169--191.

\bibitem{MartinDQC} G.J. Martin, {\em Discrete quasiconformal groups that are not the quasiconformal conjugates of M\"obius groups},  Ann. Acad. Sci. Fenn. Ser. A I Math., {\bf 11}, (1986),   179--202.


\bibitem{MLG} G.J. Martin, {\em The Hilbert-Smith conjecture for quasiconformal actions},  Electron. Res. Announc. Amer. Math. Soc., {\bf 5}, (1999), 66--70.

\bibitem{MS} G.J. Martin and R. Skora, {\em Group actions on \({\bf S}^{2} \)},  Amer. J. Math.,   {\bf 111},  (1989),  387--402.

\bibitem{OM} O.  Martio, {\em Definitions for uniform domains},  Ann. Acad. Sci. Fenn. Ser. A I Math., {\bf 5}, (1980),   197--205.

 \bibitem{MS} O. Martio and J. Sarvas, {\em Injectivity theorems in plane and space},  Ann. Acad. Sci. Fenn. Ser. A I Math., {\bf 4}, (1979),  383 -- 401.
 
\bibitem{PJ1} P. W. Jones, {\em Extension theorems for BMO},  Indiana Univ. Math. J., {\bf 29}, (1980), 1--66.

\bibitem{PJ2}  P. W. Jones, {\em  Quasiconformal mappings and extendability of functions in Sobolev spaces},  Acta Math., {\bf 147},  (1981),   71--88.

\bibitem{SK} S.P. Kerckhoff, {\em The Nielsen realization problem}, 
Ann. of Math., {\bf  117},  (1983),   235--265.

\bibitem{KM}  F. Kwakkel and V. Markovic, {\em Quasiconformal homogeneity of genus zero surfaces}, J. Anal. Math., {\bf  113},  (2011), 173--195. 
 
\bibitem{P} A. Pfluger, {\em \"Uber die \`Aquivalenz der geometrischen und der analytischen D\'efinition quasikonformer
Abbildungen},  Comment. Math. Helv., {\bf 33},  (1959), 23--33.
  

\bibitem{Tukia} P. Tukia, {\em Homeomorphic conjugates of Fuchsian groups},  J. Reine Angew. Math., {\bf 391},  (1988), 1--54.

\bibitem{Tukia2}  P. Tukia, {\em On two-dimensional quasiconformal groups}, Ann. Acad. Sci. Fenn. Ser. A I Math., {\bf 5},  (1980), 73--78.

\bibitem{Tukia3} P. Tukia, {\em  A quasiconformal group not isomorphic to a Mšbius group},  Ann. Acad. Sci. Fenn. Ser. A I Math., {\bf 6},  (1981),  149--160.
 
 \bibitem{V1} J. V\"ais\"al\"a, 
{\em On quasiconformal mappings in space},  
Ann. Acad. Sci. Fenn. Ser. A I,  (1961),  36 pp. 

 
\bibitem{V2} J. V\"ais\"al\"a, 
{\em On quasiconformal mappings of a ball}, 
Ann. Acad. Sci. Fenn. Ser. A I,  (1961), 7 pp. 
 
 \end{thebibliography}
\end{document}